\documentclass[12pt,a4paper]{amsart}
\usepackage{amssymb}
\usepackage[final]{showkeys}
\usepackage{microtype}
\usepackage{color}
\usepackage{graphicx}
\usepackage{epstopdf}

\title{The Wright functions of the second kind\\ 
    in Mathematical Physics
    \footnote{Paper published  in
MATHEMATICS (MDPI), Vol 8 No 6 (2020),  884/26pp.
  DOI: 10.3390/math8060884}
    }
\author{Francesco Mainardi$^1$}
\address{${}^1$Dipartimento di Fisica e Astronomia (DIFA), 
    	 University of  Bologna ``Alma Mater Studiorum", and INFN, \\
    	 Via Irnerio 46, I-40126 Bologna, Italy}
\email{francesco.mainardi@bo.infn.it  \  (Corresponding Author)}
    
\author{Armando Consiglio$^2$}
\address{${}^2$
Institut f\"{u}r Theoretische Physik und Astrophysik, 
\\ Universit\"{a}t W\"{u}rzburg, D-97074 W\"{u}rzburg, Germany}
\email{armando.consiglio@physik.uni-wuerzburg.de}

\def\pni{\par \noindent}

\def\vsp{\smallskip\pni}  

\def\eg{{\it e.g.}\ } \def\ie{{\it i.e.}\ }

\def\sgn{\hbox{sign}\,}

\def\e{\hbox{e}}
\def\exp{\hbox{exp}}
\def\ds{\displaystyle}
\def\RR{\mathbb{R}}
\def\Ai{\hbox{Ai}}
\def\MM{\mathbb{M}}
\def\e{{\rm e}}
\def\exp{{\rm exp}}
\def\ds{\displaystyle}

\def\q{\quad}	 \def\qq{\qquad}

\def\l{\left} \def\r{\right}

\def\rec#1{\frac{1}{#1}}

\def\L{{\mathcal L}} 
\def\F{{\mathcal F}} 
\def\NN{{\rm I\hskip-2pt N}}
\def\RR{\vbox {\hbox to 8.9pt {I\hskip-2.1pt R\hfil}}}
\def\CC{{\rm C\hskip-4.8pt \vrule height 6pt width 12000sp\hskip 5pt}}

\newcommand*{\diff}{\mathop{}\!\mathrm{d}}  

\begin{document}
	
\maketitle	
	\pagestyle{myheadings}
\markboth{\hfill \sc F. Mainardi and  A. Consiglio  \hfill}
{\hfill \sc  The  Wright functions of the second kind in Mathematical Physics \hfill}
\section{abstract}
\noindent
In this review paper we stress the importance of the higher transcendental Wright functions of the second kind in the framework of Mathematical Physics.
We first start with the analytical properties of the classical Wright functions 
of which we distinguish two kinds.
We then justify the relevance  of the Wright functions of the second kind
as fundamental solutions of the time-fractional diffusion-wave equations.
 Indeed, we think that this approach is  
the most accessible point of view.for describing 
 Non-Gaussian stochastic processes and the transition from sub-diffusion processes to wave propagation.   
  Through the sections of  the text and  suitable appendices we plan  to address
 the reader in this pathway towards the applications of the Wright functions of the second kind.\\
\\
{\it Keywords}: {Fractional Calculus, Wright Functions, Green's Functions, Diffusion-Wave Equation, Laplace Transform.}
\\ 
{\it MSC}: {26A33, 33E12, 34A08, 34C26.}

\section{Introduction}
 The special functions play a fundamental role in all  fields of Applied Mathematics and Mathematical Physics because any analytical results are  expressed in terms of 
 some of these functions. Even if the topic of special functions can appear  boring and their properties mainly treated in handbooks, we would promote the relevance
 of some of them not yet so well known. We devote our attention  to the Wright functions, in particular  with the class of the second kind.
 These functions,  as we will see hereafter,
 are fundamental to deal with some non-standard  deterministic and stochastic processes. Indeed  the Gaussian   function (known as the normal probability distribution)  must be generalized in a suitable way in the framework of partial differential equations of non-integer order.
 \\
   This work is organized as follows.
In Section 2 we introduce the Wright functions, entire in the complex plane   that we distinguish in two kinds in relation  on the value-range of the two 
 parameters on which they depend. 
 In particular we devote our attention on two Wright functions of the second kind introduced by Mainardi with the term of auxiliary functions.
 One of them, known as M-Wright function  generalizes the
 Gaussian function so it is expected to play a fundamental role in non-Gaussian stochastic processes.   
 \\
  Indeed In Section 3  we show how the Wright functions of the second kind  are relevant in the analysis of time-fractional diffusion and diffusion-wave equations
  being related to their fundamental solutions.
  This analysis leads to  generalize the known results r of the standard diffusion equation in the one-dimensional case, that is recalled in Appendix A
   by means of  auxiliary functions as particular cases of  the Wright 
   functions of the second kind known as M-Wright or Mainardi functions.
  For readers' convenience, in Appendix B   we will also  provide a  introduction to the time-derivative of fractional order in the Caputo sense
We remind that nowadays, as usual,  by fractional order we mean a non-integer order,so that the term "{\it fractional}" is a misnomer kept only for historical reasons.
\\
 In Section 4 we consider again the Mainsrdi auxiliary functions functions for their
 role in probability theory and in particular  in the framework of L{\'e}vy stable distributions whose general theory is recalled in Appendix C.
 \\
  In Section 5 we show how the  auxiliary functions turn out to be included in a class that we denote  
 {\it the four sister functions}.
 On their turn these four functions
  depending on a real parameter $\nu \in (0,1)$
are the natural generalization of 
 {\it the  three sisters functions} introduced in Appendix A devoted to the standard diffusion equation.
 The attribute of sisters was put by one of us (F. M.)  because of their inter-relations, 
 in his lecture notes on Mathematical Physics, so it has only a personal reason that we hope to be shared by the readers. 
 \\  
 Finally, in Section 6, we  provide some concluding remarks paying attention to work to be done in the next future.
 \\
 We point out that we have equipped our theoretical analysis with several plots
 hoping they will be considered illuminating for the interested readers.
 We also note that we have limited our review to the simplest boundary values problems of equations in one space dimension referring the readers to 
  suitable references for more general treatments in Section 3.1.

 \section{The Wright functions of the second kind  and the Mainardi 
auxiliary functions}  
\noindent
The classical \emph{Wright function}, that we denote by $W_{\lambda , \mu}(z)$, is defined by the series representation convergent in the whole complex plane,
\begin{equation}
W_{\lambda , \mu}(z) :=\sum_{n=0}^{\infty}{\frac{z^n}{n!\Gamma (\lambda n + \mu)}}, ~~~ \lambda > -1, ~~~ \mu \in \mathbb{C},
\end{equation}
The \emph{integral representation} reads as: 
\begin{equation}
W_{\lambda , \mu}(z) = \frac{1}{2\pi i}\int_{Ha_{-} }{e^{\sigma + z\sigma ^{-\lambda}}\frac{d\sigma}{\sigma ^{\mu}}}, ~~~ \lambda > -1, ~~~ \mu \in \mathbb{C},
\end{equation}
where $Ha_{-}$ denotes the Hankel path: this one is a loop which starts from $-\infty$ along the lower side of negative real axis, encircles with a small circle the axes origin and ends at $-\infty$ along the upper side of the negative real axis.
\\
$W_{\lambda , \mu}(z)$ is then an \emph{entire function} for all
$\lambda \in (-1, +\infty)$.
 Originally Wright assumed $\lambda \ge 0$ in connection with his investigations on the asymptotic theory of partition
 \cite{Wright 33,Wright 35} and only in 1940 he considered
  $-1 < \lambda < 0$, \cite{Wright 40}.
We note that in the Vol 3, Chapter 18 of the handbook of the Bateman Project
\cite{Erdelyi BATEMAN}, presumably for a misprint,
 the parameter $\lambda$  is restricted to be non-negative,
whereas  the Wright functions remained  practically ignored   in other  handbooks.  
In 1993  Mainardi,   being aware only  of  the Bateman handbook,
proved that the Wright function is entire also for $-1<\lambda<0$  in his approaches to the time fractional diffusion equation, that will be dealt in a next Section.
\\
In view of the asymptotic representation in the complex domain 
and of the Laplace transform for positive argument  $z=r>0$
($r$ can be the time variable $t$ or the space variable $x$)
  the Wright functions are distinguished in
   \emph{first kind} ($\lambda \geq 0$) and \emph{second kind}
 ($-1< \lambda < 0$) 
 as outlined in the Appendix F of the book by Mainardi
 \cite{Mainardi BOOK2010}.
 In particular,  for the asymptotic behaviour, we refer the interested reader to
 the two  papers by Wong and Zhao \cite{Wong 99a,Wong 99b},
  and to  the surveys by Luchko and by Paris in the Handbook of Fractional Calculus
  and Applications,
   see respectively \cite{Luchko HFCA}, \cite{Paris HFCA},
    and references therein.
   \\ 
We note that the Wright functions  are  entire of order $1/(1+\lambda)$
 hence only  the first kind  functions ($\lambda \ge 0$)  are of exponential order whereas  the second kind functions ($-1<\lambda<0$)  
 are not of exponential\ order.
 The case $\lambda =0$ is trivial since
 $W_{0,\mu}(z) = {\e^z}/{\Gamma(\mu)}.$
As a consequence of the distinction in the kinds, we must point out the different Laplace transforms proved 
e.g. in  \cite{GOLUMA 99},\cite{Mainardi BOOK2010}, see also the 
recent survey
on Wright functions by Luchko \cite{Luchko HFCA}. 
We have:
\\
for the first kind, when $\lambda \ge 0$
\begin{equation}
  W_{\lambda ,\mu } (\pm r) \,\div \,
    \rec{s}\, E_{\lambda ,\mu }\left(\pm \rec{s}\right) \,;
\end{equation}
\\
for the second kind,   when $-1<\lambda<0$
and putting for convenience  $\nu = -\lambda$
so $0< \nu<1$
\begin{equation}
W_{-\nu  ,\mu } (-r) \,\div \,
    E_{\nu ,\mu+\nu }\left(-s \right) \,.
\end{equation}
Above we have  introduced   
the  Mittag-Leffler function
in two parameters $\alpha>0$, $\beta \in \CC$ 
  defined as its convergent series for  all  $z \in \CC$
  \begin{equation}
E_{\alpha, \beta}(z) := \sum_{n=0}^{\infty}
{\frac{z^n}{\Gamma(\alpha n + \beta)}}.
\label{eq:mittag-leffler}
\end{equation}
For more details on the special functions of the Mittag-leffler type we refer the interested readers to the treatise by Gorenflo et al 
\cite{GKMS BOOK14},
where in the forthcoming 2-nd edition also the Wright functions are treated 
in some detail. 
\\ 
In particular, two Wright functions of the second kind, 
originally introduced by Mainardi and 
named $F_\nu (z)$ and $M_\nu (z)$ ($0<\nu<1$),  are called 
\emph{auxiliary functions} in virtue of their role in the time fractional  diffusion equations considered in the next section.
 These functions, $F_\nu (z)$ and $M_\nu (z)$, are  indeed special cases of the Wright function of the second kind $W_{\lambda , \mu}(z)$ by setting, respectively,
  $\lambda = -\nu$ and $\mu = 0$ or $\mu = 1-\nu$. 
  Hence we have:
\begin{equation}
F_{\nu}(z) := W_{-\nu , 0}(-z), ~~~~ 0 < \nu < 1,
\end{equation}
and
\begin{equation}
M_{\nu}(z) := W_{-\nu , 1- \nu}(-z),~~ 0 < \nu < 1,
\end{equation}
Those functions are interrelated through the following relation: 
\begin{equation}
F_{\nu}(z) = \nu z M_{\nu}(z),
\end{equation}
which reminds us the second relation in \ref{eq:reciprocityrelationstandarddiff}, seen for the standard diffusion equation.
\\
The series representations of the auxiliary functions are derived from those of $W_{\lambda , \mu}(z)$. Then:
\begin{equation}
F_{\nu}(z):=  \sum_{n=1}^{\infty}{\frac{(-z)^n}{n!\Gamma(-\nu n)}}
= - \frac{1}{\pi}\sum_{n=1}^{\infty}{\frac{(-z)^{n}}{n!}\Gamma(\nu n + 1)\sin{(\pi \nu n)}},
\label{eq:F}
\end{equation}
and
\begin{equation}
M_{\nu}(z) :=  \sum_{n=0}^{\infty}{\frac{(-z)^n}{n!\Gamma[-\nu n + (1-\nu )]  }} 
= \frac{1}{\pi}\sum_{n=1}^{\infty}{\frac{(-z)^{n-1}}{(n-1)!}\Gamma(\nu n)\sin{(\pi \nu n)}},
\label{eq:M}
\end{equation}
where it has been used in both cases the \emph{reflection formula} for the Gamma function (Eq. \ref{eq:reflectionformula}) among the first and the second step of Eqs. (\ref{eq:F}) and (\ref{eq:M}),
\begin{equation}
\Gamma(\zeta)\Gamma(1-\zeta ) = \pi / \sin{\pi \zeta}.
\label{eq:reflectionformula}
\end{equation}
{\bf Remark} :
In the present  version  we have corrected an error occurring in  version V1
on the RHS of Eq. (9) referring to the $F_\nu$ function.
Unfortunately this error was already present in the  Appendix F of the first edition 
of the 2010 book by Mainardi   \cite{Mainardi BOOK2010}.
\newpage
Also the integral representations of the auxiliary functions are derived from those of $W_{\lambda , \mu}(z)$. Then: 
\begin{equation}
F_{\nu}(z) := \frac{1}{2\pi i}\int_{Ha_{-}}{e^{\sigma - z\sigma ^{\nu}}d\sigma}, ~~~ z \in \mathbb{C}, ~~~ 0 < \nu < 1 
\end{equation}
and
\begin{equation}
M_{\nu}(z) := \frac{1}{2\pi i}\int_{Ha_{-}}{e^{\sigma - z\sigma ^{\nu}}\frac{d\sigma}{\sigma ^{1-\nu}}}, ~~~ z \in \mathbb{C}, ~~~ 0 < \nu < 1 
\end{equation}
\medskip
Explicit expressions of $F_{\nu}(z)$ and $M_{\nu}(z)$ in terms of known functions are expected for some particular values of $\nu$
as shown and recalled by Mainardi in the first 1990's  in a series of papers 
\cite{Mainardi WASCOM93, Mainardi RADIOPHYSICA95,Mainardi AML96,Mainardi CHAOS96}, that is
\begin{equation}
M_{1/2}(z) = \frac{1}{\sqrt{\pi}}e^{-z^2 /4},
\end{equation}
\begin{equation}
M_{1/3}(z) = 3^{2/3}\Ai(z/3^{1/3}),. 
\end{equation}
Liemert and Klenie \cite{Liemert-Kleine JMP2015} have added the 
following expression for $\nu=2/3$
\begin{equation}
M_{2/3}(z) = 3^{-2/3}
\left[3^{1/3}\,z\, \Ai\left(z^2/3^{4/3}\right) -
3\Ai^\prime\left ( z^2/3^{4/3}\right)\right ]
\, \e^{-2z^3/27},
\end{equation}
where $\Ai$  and $\Ai^\prime$ denote the \emph{Airy function}
and its first derivative.
Furthermore they  have 
suggested in the positive real field  $\RR^+$ the following remarkably integral representation
\begin{equation}
M_\nu(x) =
\frac{1}{\pi}\,
\frac{x^{\nu/(1-\nu)}}{1-\nu}\,
\int_0^\pi \! C_\nu(\phi)\, \exp \left(-C_\nu(\phi)\right)\,
x^{1/(1-\nu)}\, d\phi,
\end{equation}
where
\begin{equation}
C_\nu(\phi) = \frac{\sin(1-\nu)}{\sin \phi}\,
\left( \frac{\sin \nu \phi}{\sin \phi}\right)^{\nu/(1-\nu)}\,,
\end{equation}
corresponding to equation (7) of the article written by Saa and Venegeroles
\cite{Saa-Venegeroles PRE2011} .
\\
Let us point out     the asymptotic behaviour 
of the function $M_\nu(x)$ as $x \to +\infty$. 
Choosing as
a variable $x/\nu $ rather than $x$, the computation of the 
asymptotic representation  by the saddle-point approximation
 carried out by Mainardi and Tomirotti yields,  
 see \cite{Mainardi-Tomirotti TMSF94} and 
 \cite{Mainardi BOOK2010},
\begin{equation} 
\! M_\nu (x/\nu ) \sim
   a(\nu )\, x^{\ds{(\nu -1/2)/(1-\nu)}}
  \,
   \exp \left[- b(\nu)\,x^{\ds {1/(1-\nu)}}\right],
\end{equation}
where
\begin{equation}
  a(\nu) = \rec{\sqrt{2\pi\,(1-\nu)}}    >0 \,,  \q
  b(\nu) = \frac{1-\nu }{  \nu }    >0 \,. 
  \end{equation}
The above evaluation is consistent with the first term in
the original  asymptotic expansion  by Wright in 
\cite{Wright 35, Wright 40} after having
used the  definition of .$M$-Wright function 
 \\
Now we find it convenient to show the plots of the $M$-Wright functions 
on a space symmetric interval of $\RR$  in Figs 1, 2, corresponding to the cases 
$0\leq \nu \le 1/2$ and 
$1/2 \leq \nu \le 1$, respectively.
   We recognize 
the non-negativity of the $M$-Wright function on $\RR$ for
$1/2 \leq \nu \leq 1$  consistently with the  analysis
 on distribution of zeros and asymptotics of Wright functions
 carried out by Luchko, see \cite{Luchko 2000}, \cite{Luchko HFCA}.  
 \begin{figure}[h!]
	\centering
	\includegraphics[width=120mm, height=80mm]{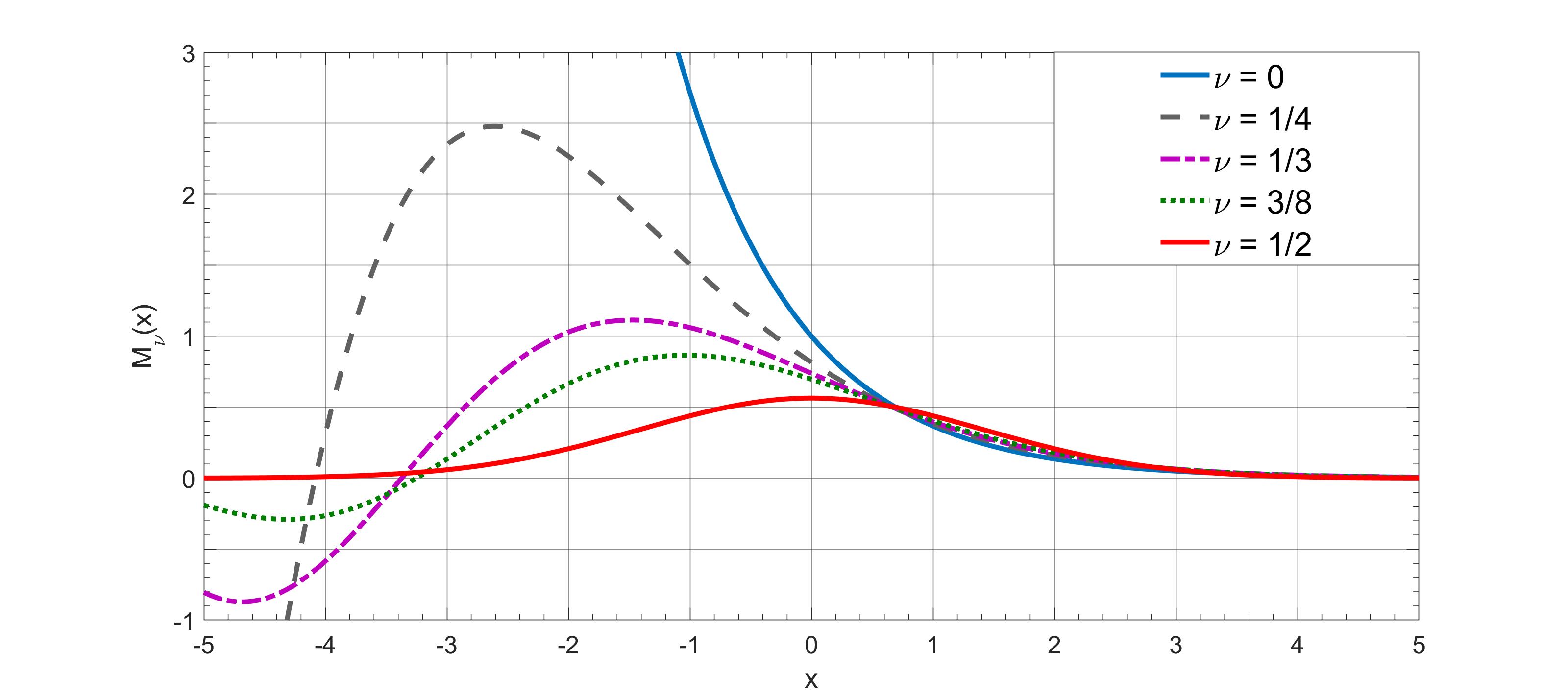}
	\vskip-0.5truecm
	\caption{Plots of the $M$-Wright function as a function of the $x$ variable, for $0 \leq \nu \leq 1/2$.}
	\label{fig:subdiff}
\end{figure}
\begin{figure}[h!]
\centering
\includegraphics[width=120mm, height=80mm]{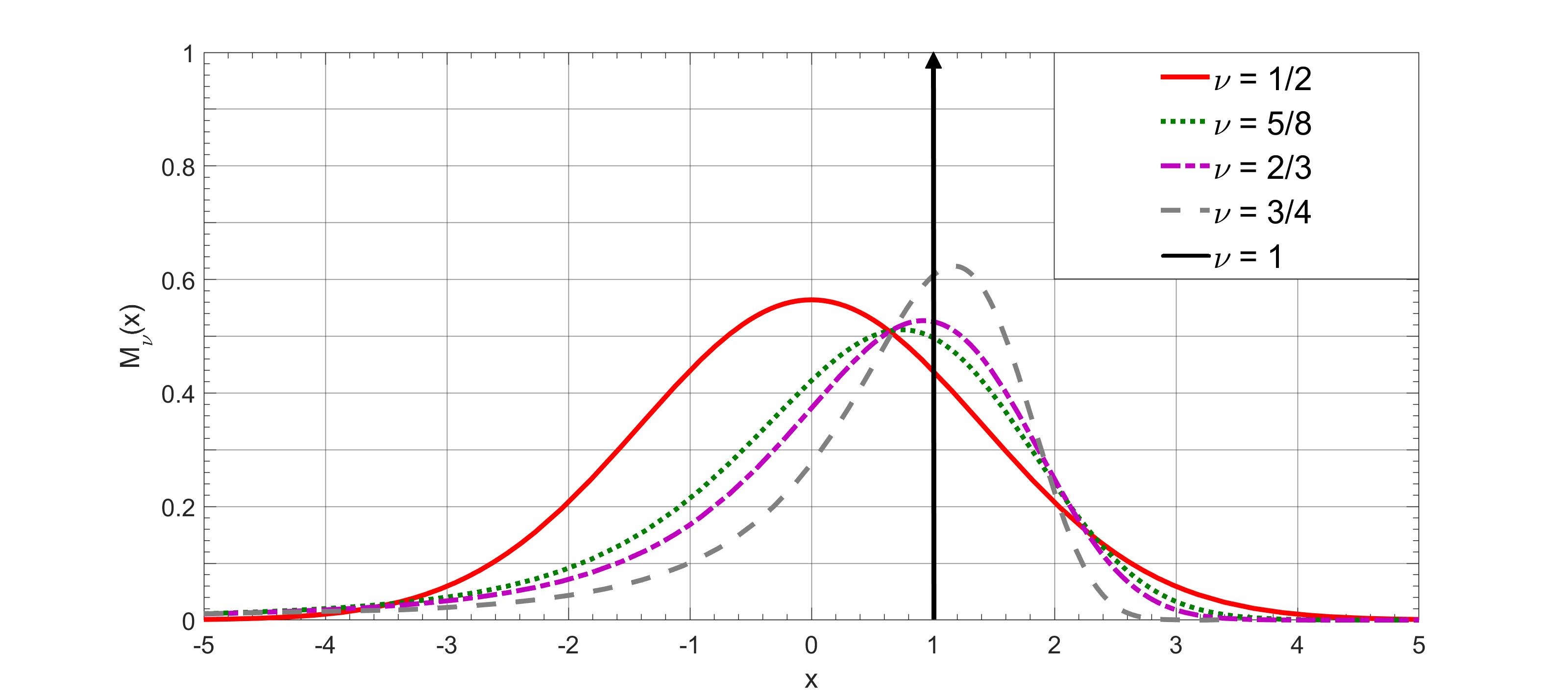}
\vskip-0.5truecm
\caption{Plots of the $M$-Wright function as a function of the $x$ variable, for $1/2 \leq \nu \leq 1$.}
\label{fig:diffwavenomodule}
\end{figure}
\section{The Wright functions of the second kind 
and the time-fractional diffusion wave equation}
\noindent
As we will see the Wright functions of the second kind are relevant in the analysis of the Time-Fractional Diffusion-Wave  Equation (TFDWE).
\\
For this purpose we introduce now the TFDWE as a generalization of the standard diffusion equation and we see how the two Mainardi auxiliary functions 
 come into play.
\\ 
The TFDWE 
is so obtained from the standard diffusion equation (or the D'Alembert wave equation) by replacing the first-order (or the second-order) time derivative by a fractional derivative (of order $0 < \beta \leq 2$) in the Caputo sense, obtaining the following Fractional PDE:
\begin{equation}
	\frac{\partial ^{\beta}u}{\partial t^{\beta}} = D\frac{\partial ^{2}u}{\partial x^{2}} ~~~~~~~~ 0 < \beta \leq 2, ~~~~ D>0,
	\label{eq:tfde}
\end{equation}
where $D$ is a positive constant whose dimensions are $L^2T^{-\beta}$ and $u=u(x, t; \beta)$ is the field variable, which is assumed again to be a causal function of time.
The Caputo fractional derivative is recalled in the Appendix B so that in explicit form the TFDWE  (\ref{eq:tfde}) 
splits in the following integro-differential equations:
\begin{equation}
\frac{1}{\Gamma (1-\beta)}\int_{0}^{t}{(t-\tau)^{-\beta}\biggl( \frac{\partial u}{\partial \tau}\biggl)d\tau} = D\, \frac{\partial ^2 u}{\partial x^2}, ~~~ 0 < \beta \leq 1;
\end{equation}
\begin{equation}
\frac{1}{\Gamma (2-\beta)}\int_{0}^{t}{(t-\tau)^{1-\beta}\biggl( \frac{\partial^2  u}{\partial \tau ^2}\biggl)d\tau} = D\, \frac{\partial ^2 u}{\partial x^2}, ~~~ 1 < \beta \leq 2.
\end{equation}
\\
In view of our analysis we find convenient to put: 
\begin{equation}
\nu = \frac{\beta}{2}, ~~~ 0 < \nu \leq 1.
\end{equation}
\\
We can then formulate the basic problems for the Time Fractional Diffusion-Wave Equation using a correspondence with the two problems for the standard diffusion equation.
\\
Denoting by $f(x)$ and $g(t)$ two given, sufficiently well-behaved functions, we define:
\\
a) Cauchy problem 
\begin{equation}
\begin{cases} u(x, 0^+; \nu) = f(x), ~~~~~~~~~~~~~~~~ -\infty < x < +\infty; & \mbox{} \mbox{ } \\ u(\pm\infty , t;\nu) = 0,  ~~~~~~~~~~~~~~~~~~~~~~ t > 0 & \mbox{} \mbox{} \end{cases}
\end{equation}
\medskip
b) Signalling problem
\begin{equation}
\begin{cases} u(x, 0^+; \nu) = 0, ~~~~~~~~~~~~~~~~~~~~~~~~ 0 \leq x < +\infty; & \mbox{} \mbox{ } \\ u(0^+ , t;\nu) = g(t), ~u(+\infty , t;\nu) = 0,  ~~ t > 0 & \mbox{} \mbox{} \end{cases}
\end{equation}
\\
If  $1/2 < \nu \leq 1$ corresponding to $1< \beta \leq 2$ 
we must consider also the initial value of the first time derivative of the field variable $u_t (x, 0^+ ; \nu)$, since in this case Eq. (\ref{eq:tfde})  turns out to be akin to the wave equation and consequently two linear independent solutions are to be determined. However, to ensure the continuous dependence of the solutions to our basic problems on the parameter $\nu$ in the transition from 
$\nu = (1/2)^-$ to $\nu = (1/2)^+$, we agree to assume $u_t (x, 0^+ ; \nu) = 0$.\\
\\
For the Cauchy and Signalling problems,
following the approaches by Mainardi, see e.g. \cite{Mainardi WASCOM93} and related papers,
we introduce now the Green functions 
$\mathcal{G}_c (x, t;\nu )$ and $\mathcal{G}_s (x, t;\nu )$ that for both problems can be determined by the $LT$ technique, so extending the results known from the ordinary diffusion equation.
We recall that the Green functions are also referred to as the fundamental solutions, corresponding  respectively to $f(x) =\delta(x)$ and $g(t) = \delta(t)$
with $\delta(\cdot)$ is the Dirac delta generalized function 
\\
The expressions for the Laplace Transforms of the two Green's functions are:
\begin{equation}
\widetilde{\mathcal{G}} _c (x, s; \nu) = \frac{1}{2\sqrt{D}s^{1-\nu}}
\, \e^{(-|x|/\sqrt{D})s^{\nu}}
\label{eq:laplacetransformcauchygreenfunctionfractionalcalculus}
\end{equation}
and
\begin{equation}
\widetilde{\mathcal{G}}_s (x, s ; \nu ) = \e^{-(x/\sqrt{D})s^{\nu}}
\end{equation}
\\
Now we can easily recognize the following relation:
\begin{equation}
\frac{d}{ds}\widetilde{\mathcal{G}}_s = 
-2\nu \, x\, \widetilde{\mathcal{G}}_c , ~~~ x > 0
\end{equation}
which implies for the original Green functions the following \emph{reciprocity relation} for $x>0$`and $t>0$ and $0<\nu<1$:
\begin{equation}
	2\nu x \mathcal{G}_c (x, t ; \nu ) = t \mathcal{G}_s (x, t ; \nu ) = F_{\nu}(z) = \nu z M_{\nu} (z)\,, \quad 
	z = \frac{x}{\sqrt{D}t^{\nu}}
	\label{eq:reciprocityrelationfractionalcalculus}
\end{equation}  
where $z$ is
 the \emph{similarity variable}
 and $F_\nu(z)$ and $M_\nu(z)$ are the Mainardi auxilary functions
 introduced in the previous section.
Indeed Eq. (\ref{eq:reciprocityrelationfractionalcalculus}) is the
generalization of Eq. (A.8) that we have seen for the standard diffusion equation
due to the introduction of the time fractional derivative of order $\nu$
\\
Then,  the two Green functions  of  the Cauchy and Signalling problems
turn out to be expressed in terms of the two auxiliary functions as follows.
\\
 For the Cauchy problem we have
\begin{equation}
	\mathcal{G}_C (x, t;\nu) = \frac{t^{-\nu}}{2\sqrt{D}}\,
	M_\nu\left ( \frac{|x|}{\sqrt{D}t^{\nu}}\right)\,, \quad -\infty<x<+\infty\,, \quad
	t \ge 0\,,
	\label{eq:Green-Cauchy-nu}
\end{equation}
that generalizes Eq. (A.5). 
\\
For the Signalling problem we have:
\begin{equation}
	\mathcal{G}_S (x, t;\nu) = \frac{\nu x t^{-\nu-1}}{\sqrt{D}}\, 
	M_\nu\left ( \frac{x}{\sqrt{D}t^{\nu}}\right)\,, 
\quad x\ge 0, \quad t\ge 0\,,
\label{eq:Green-Signalling-nu}
\end{equation}
that generalizes Eq. (A.7).  
\subsection{Complements to the time-fractional diffusion-wave equations} 

The boundary value problems dealt previously can be considered with a source data function $fx)$ and $g(t)$ different from the Dirac generalized functions, in particular with  box-type functions as it has been carried out recently by us, 
see \cite{Consiglio-Mainardi RIM2019}.
  
The TFDWE can be generalized in 2D and 3D space dimensions. so consequently the Wright functions play again a fundamental role. However,  we prefer to refer the interested reader to the literature, in particular to the papers
 by Luchko and collaborators
 \cite{Luchko 2000,Luchko CAIM2014,Luchko MATHS2017,Luchko HFCA},
 \cite{Luchko-Mainardi JVA2014,Luchko-Mainardi HFCA},
 \cite{Boyadjiev-Luchko CSF2017},
by Hanyga \cite{Hanyga PRSA2002}
 and to  the recent analysis by Kemppainen \cite{Kemppainen 2019}.
All of them are originated in some way from the seminal paper by
 Schneider \& Wyss  \cite{Schneider-Wyss JMP1989}.
In some of  these papers the authors have considered  also fractional differentiation 
both in time and in space, so that they have  generalized to more than one dimension the former analysis  by Mainardi, Luchko \& Pagnini
\cite{Mainardi LUMAPA01}
on the space-time fractional diffusion-wave equations.    
 
\section{The $M-$Wright functions in probability theory and the stable distributions}
We  recognize that the Wright $M$-function with support  in $\RR^+$
can be interpreted as probability density function ($pdf$)
because it is non negative and also it
satisfies the normalization condition:
\begin{equation}
\int_{0}^{\infty}{M_{\nu}(x)\, dx} = 1\,.
\end{equation}
\\
  We now  provide more details on these densities in the framework of the theory of probability. 
  \\ 
Fundamental quantities about the Wright $M-$function are the \emph{absolute moments} of order $\delta > -1$  
in $\mathbb{R}^+$, that are finite and turn out to be:
\begin{equation}
\int_{0}^{\infty}{x^{\delta}\,M_{\nu}(x)\, dx} 
= \frac{\Gamma(\delta + 1)}{\Gamma(\nu \delta + 1)}, ~~~ \delta > -1, ~~~ 0 \leq \nu < 1.
\end{equation}
The result is based on the integral representation of the $M-$Wright function:
\\
\begin{equation}
\begin{split}
\int_{0}^{\infty}{x^{\delta}M_{\nu}(x)dx} &= \int_{0}^{\infty}{x^{\delta}\Biggl[\frac{1}{2\pi i}\int_{Ha_{-}}{e^{\sigma - x\sigma ^{\nu}}\frac{d\sigma}{\sigma ^{1-\nu}}} \Biggl]dx}\\
&=\frac{1}{2\pi i}\int_{Ha_{-}}{e^{\sigma}\Biggl[ \int_{0}^{\infty}{e^{-x\sigma ^{\nu}}x^{\delta}dx} \Biggl]\frac{d\sigma}{\sigma ^{1-\nu}}}\\
&= \frac{\Gamma(\delta + 1)}{2\pi i}\int_{Ha_{-}}{\frac{e^{\sigma}}{\sigma ^{\nu \delta + 1}}d\sigma} = \frac{\Gamma(\delta + 1)}{\Gamma(\nu \delta + 1)}
\end{split}
\label{eq:absolutemoments}
\end{equation}
The exchange between two integrals and the following identity contributed to the final result for Eq. (\ref{eq:absolutemoments}):
\begin{equation}
\int_{0}^{\infty}{\e^{-x\sigma ^{\nu}}x^{\delta}dx} = \frac{\Gamma(\delta + 1)}{(\sigma ^{\nu})^{\delta + 1}}.
\end{equation}
\\
In particular, for $\delta = n \in \NN$, the above formula provides the moments of integer order. 
Indeed
 recalling the Mittag-Leffler function introduced in
  Eq. (\ref{eq:mittag-leffler})
 with $\alpha=\nu$ and $\beta=1$:
 \begin{equation}
E_{\nu}(z) := \sum_{n=0}^{\infty}{\frac{z^n}{\Gamma(\nu n + 1)}}, ~~~ \nu > 0, ~~~ z \in \mathbb{C}, 
\end{equation}
the moments of integer order can also be computed   from the Laplace transform pair 
\begin{equation}
M_\nu(x) \div E_\nu(s)\,,
\end{equation}
proved in the Appendix F of \cite{Mainardi BOOK2010}
 as follows:
\begin{equation} 
    \int_0^{+\infty} \!\!\! x^{\, n}\, M_\nu (x)\, dx = 
   \lim_{s\to 0} \, (-1)^n \,
 \frac{d^n}{  ds^n}\, E_\nu (-s)
   = \frac{\Gamma(n+1)}{  \Gamma(\nu n+1)}\,.
   \end{equation}
\subsection{The auxiliary functions versus extremal  stable densities}
We find it 
worthwhile to recall   the relations between the Mainardi auxilary functions   and the extremal  L\'evy stable densities as proven
 in the 1997 paper by Mainardi and Tomirotti \cite{Mainardi-Tomirotti GEO97}.
For readers' convenience we  refer  to Appendix C for an essential  account of the general L{\'e}vy stable distributions in probability.
 Indeed, from a comparison between  the series expansions of stable densities in (C.8)-(C.9) and  of the  auxiliary functions in 
 Eqs. (\ref{eq:F}) - (\ref{eq:M}), 
we recognize that 
the 
auxiliary functions  are related to the extremal stable densities
as follows
\begin{equation}
L_\alpha ^{-\alpha } (x)
 =  \rec{x}\,  F_\alpha  (x^{-\alpha }) =
 {\alpha  \over x^{\alpha  +1}}\,  M_\alpha (x^{-\alpha }) \,,\q
 0<\alpha  <1\,, \q x \ge 0\,,
 \label{eq:unilateral-extremal}
  \end{equation} 
\begin{equation}  
  L_\alpha ^{\alpha  -2}(x)
 = \rec{x}\,  F_{1/\alpha }(x) =
 {1 \over \alpha  }\,  M_{1/\alpha }(x) \,,\q
 1<\alpha  \le 2\,, \q -\infty<x<+\infty\,
 \label{eq:bilateral-extremal}.
 \end{equation}
 In the above equations, for $\alpha=1$, the skewness parameter turns out to be
  $\theta =-1$, 
 so  we get the singular limit
 \begin{equation} 
 L_1^{-1}(x)= M_1(x)= \delta(x-1)\,.
 \end{equation}
Hereafter we show the plots the extremal stable densities according to their expressions in terms of the $M$-Wright functions, see   
Eq. (\ref{eq:unilateral-extremal})., Eq. (\ref{eq:bilateral-extremal})
for $\alpha=1/2$ and
 $\alpha=3/2$, respectively.
 \begin{figure}[h!]
	\centering
	\includegraphics[width= 120mm, height=50mm]{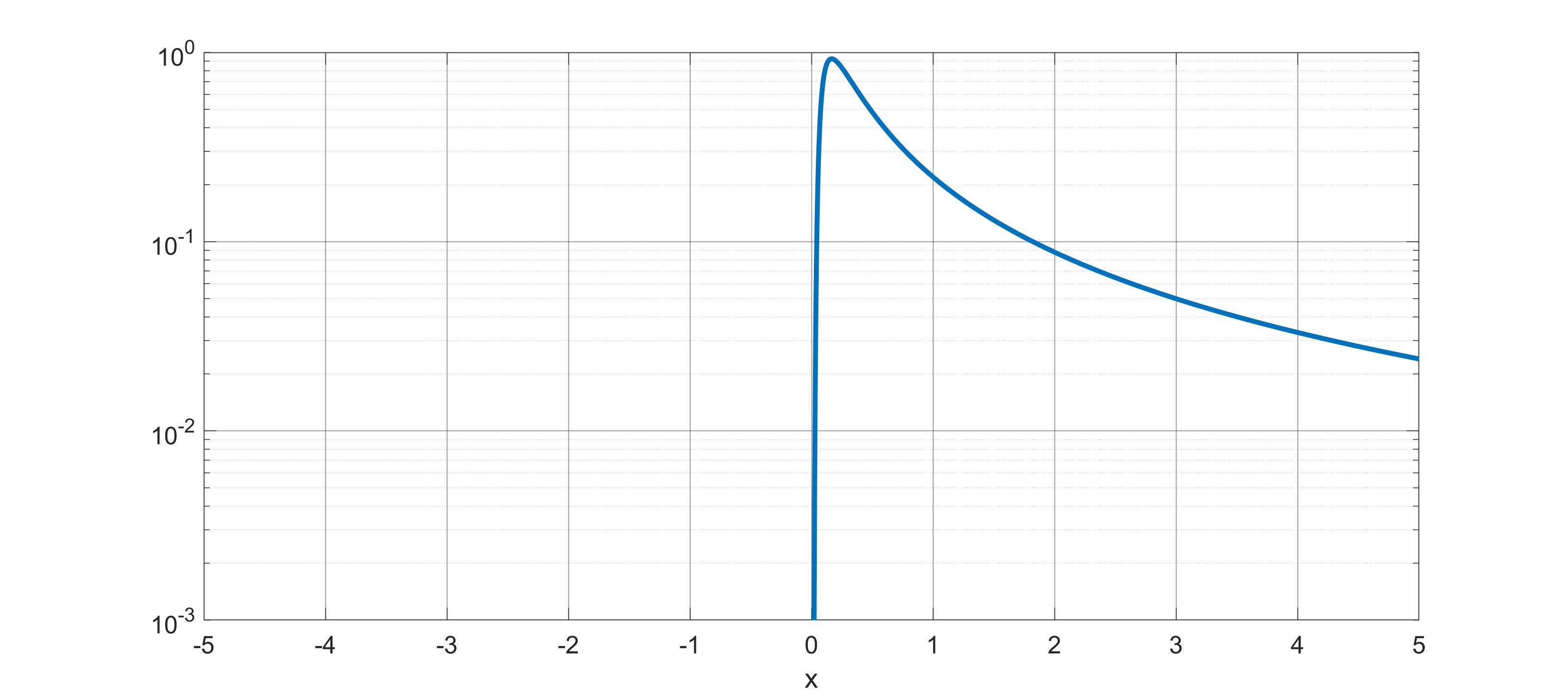}
\vskip -0.2truecm
	\caption{Plot of the unilateral extremal stable pdf for $\alpha=1/2$} 
	 \end{figure}
	 \begin{figure}[h!]
	\centering
	\includegraphics[width= 120mm, height=50mm]{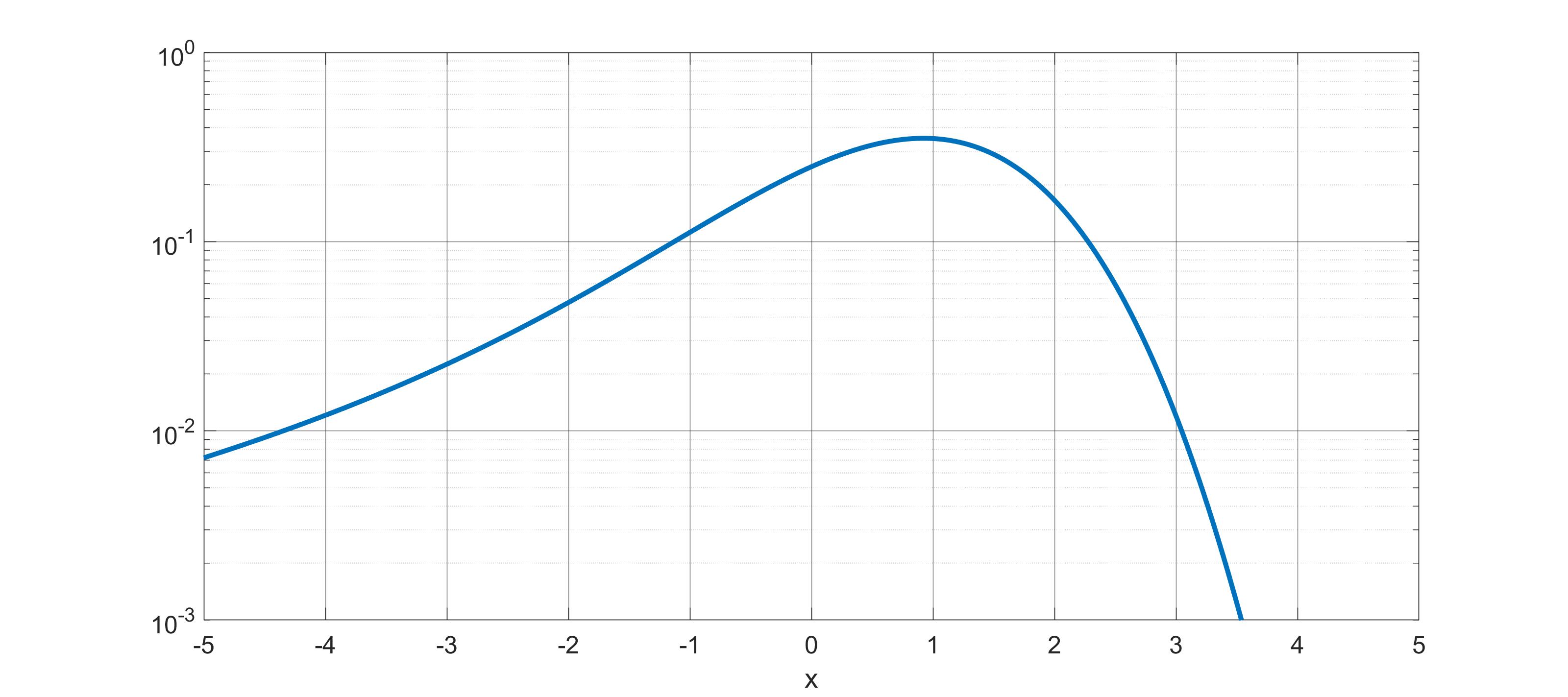}
	\vskip -0.2truecm
	\caption{Plot of  the bilateral  extremal stable pdf for $\alpha=3/2$} 
	 \end{figure}
	 We recognize that the above plots are consistent with the corresponding 
	 ones shown by Mainardi et al. \cite{Mainardi LUMAPA01} 
	 for the stable pdf's derived 
	 as fundamental solutions of  a suitable space-fractional diffusion equation. 
\subsection{The symmetric $M-$ Wright function} 
We easily recognize that 
extending the function $M_\nu(x)$
in a symmetric way to all of $\RR$ 
(that is putting $x=|x|$) and dividing by 2
 we have  a {\it symmetric} $pdf$  
 with support in all of $\RR$.
\\
 As the parameter $\nu$ changes between 0 and 1, 
the {\it pdf} goes from the Laplace {\it pdf} 
 to two half discrete delta {\it pdf}s  passing for $\nu=1/2$
  through the Gaussian {\it pdf}. 
     \\
To develop a visual intuition, also in view of the subsequent applications, 
we show the plots of the symmetric $M-$Wright function on the real axis at $t=1$ for some rational values of the parameter $\nu \in [0, 1]$
\begin{figure}[h!]
	\centering
	\includegraphics[width= 120mm, height=70mm]{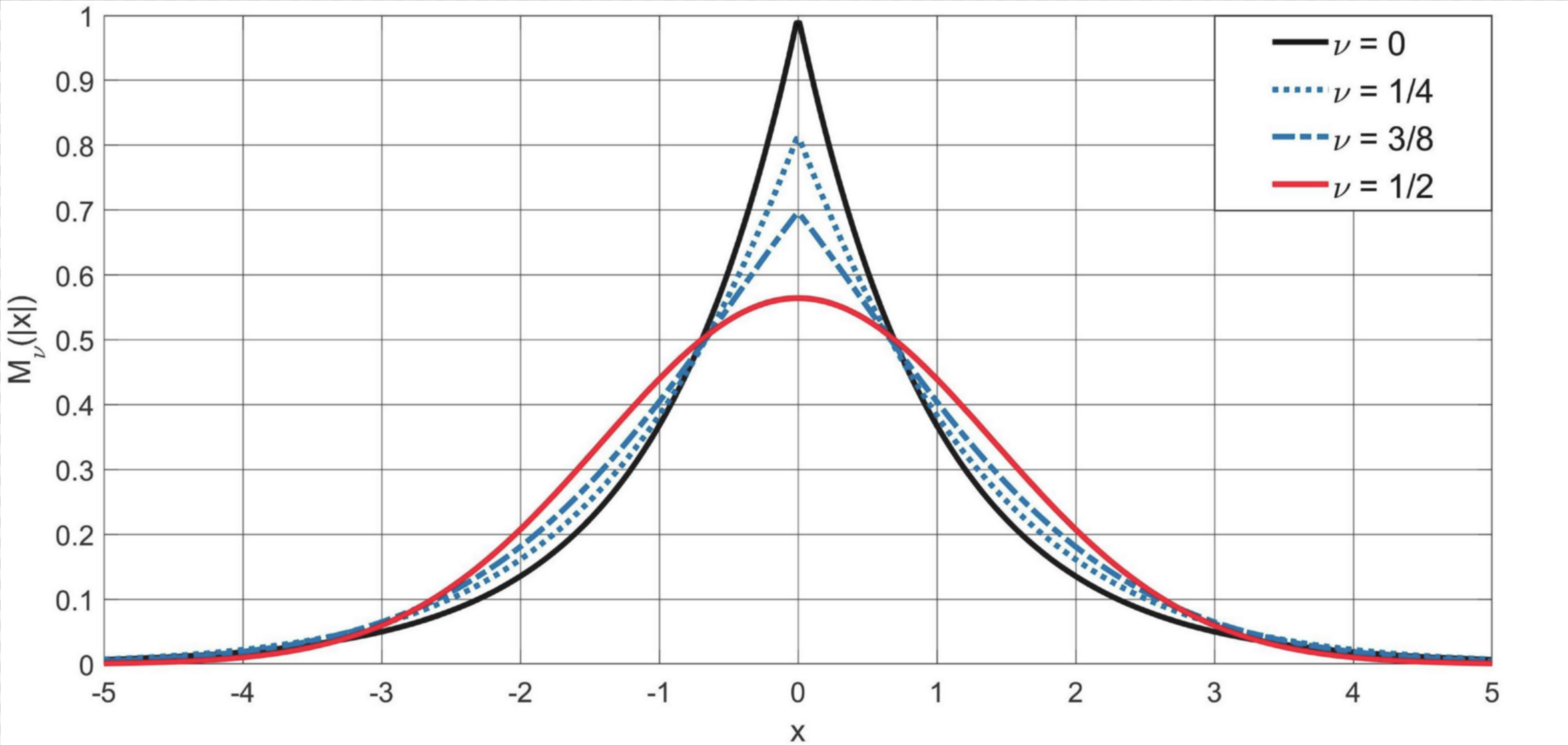}
	\caption{Plot of the symmetric $M-$Wright  function $M_{\nu}(|x|)$ 
	for $0 \leq \nu \leq 1/2$. Note that the $M-$Wright function becomes
a Gaussian density  for $\nu=1/2$.}	
	 \end{figure}
\begin{figure}[h!]
	\centering
	\includegraphics[width= 120mm, height=70mm]{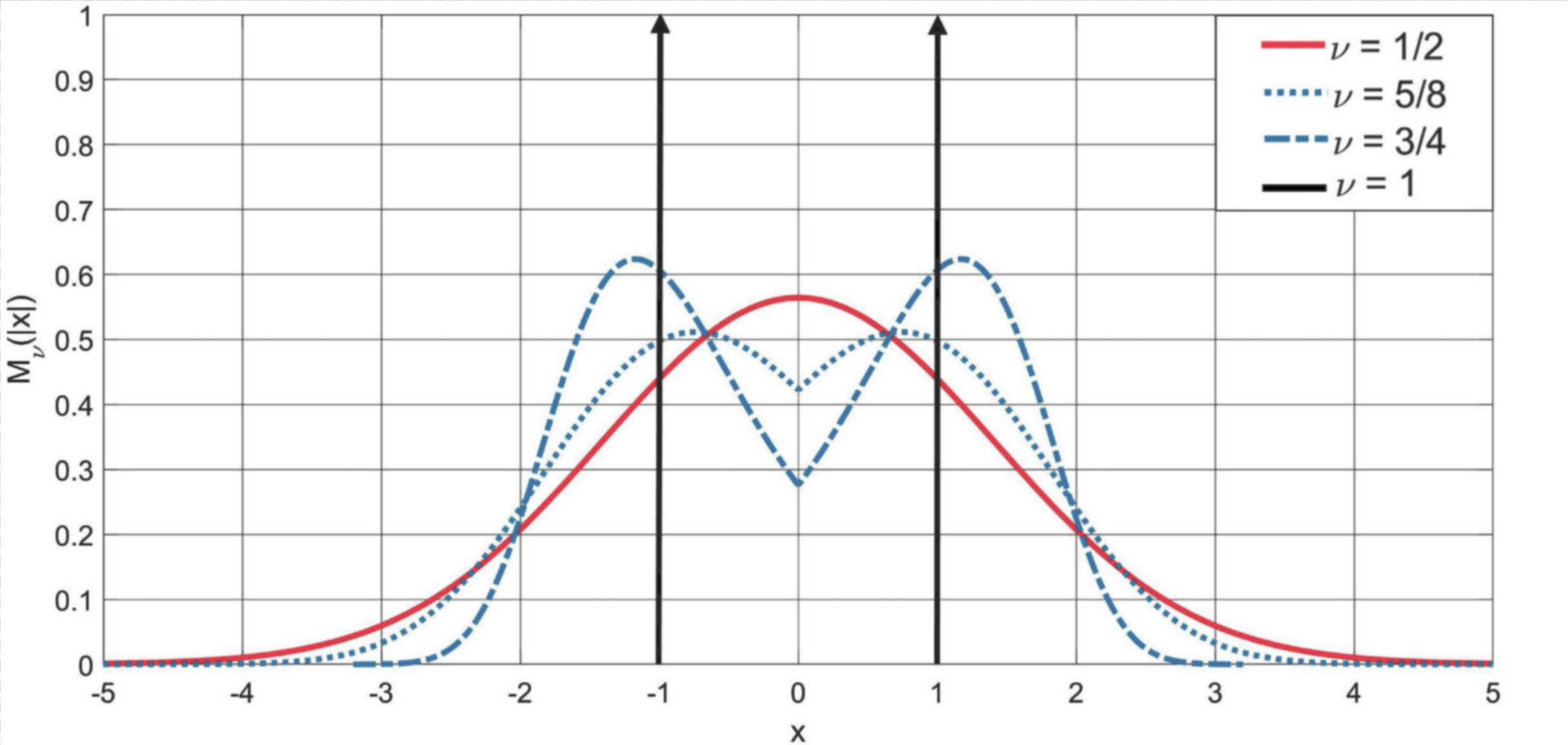}
	\caption{Plot of the symmetric $M-$Wright type function $M_{\nu}(|x|)|$
	 for $ 1/2 \leq \nu \leq 1$. Note that the $M$Wright function becomes a
a sum of two delta functions centered in $x=\pm 1$ for $\nu = 1$.}	 
\end{figure}
\\
The  readers are invited to look the  YouTube  video 
by Consiglio whose title is ``Simulation of the  $M-$Wright function", in which 
the author   shows  the evolution of this  function as the parameter $\nu$ changes between 0 and 0.85 in a finite  interval of $\RR$ centered in $x=0$.
\\
Finally we compute the \emph{characteristic function} for the symmetric 
$M-$Wright 
{\it pdf}.  We get for $M_\nu(|x|$ 
\begin{equation}
\begin{split}
\mathcal{F}\Bigl[\frac{1}{2}M_{\nu}(|x|)\Bigl] &:= \frac{1}{2}\int_{-\infty}^{+\infty}{e^{+i\kappa x}M_{\nu}(|x|)dx}\\
& ~ =\int_{0}^{\infty}{\cos{(\kappa x)}M_{\nu}(x)dx} = E_{2\nu}(-\kappa ^2)
\end{split}
\label{eq:characteristicfunction}
\end{equation}

\medskip
Eq. \ref{eq:characteristicfunction} is obtained by developing in series the cosine function and using Eq. \ref{eq:absolutemoments}. In particular:
\begin{equation}
\begin{split}
\int_{0}^{\infty}{\cos{(\kappa x)}M_{\nu}(x)dx} &= \sum_{n=0}^{\infty}{(-1)^n \frac{\kappa ^{2n}}{(2n)!}\int_{0}^{\infty}{x^{2n}}M_{\nu}(x)dx}\\
&=\sum_{n=0}^{\infty}{(-1)^n\frac{\kappa ^{2n}}{\Gamma (2\nu n + 1)}} = E_{2\nu}(-\kappa ^2)
\end{split}
\end{equation}


 \subsection{The Wright $\MM$-function in two variables.}
In view of time-fractional diffusion processes related to time-fractional diffusion equations 
it is worthwhile to introduce the 
function  in two variables
\begin{equation} 
  \MM_\nu(x,t):= t^{-\nu}\, M_\nu(xt^{-\nu})\,,\q 0<\nu < 1\,,\q x,t \in \RR^+ \,,
  \end{equation}   
  which defines a  spatial probability density in $x$ evolving in 
  time $t$ with self-similarity exponent $H=\nu$.
  Of course for $x\in \RR$ 
  we have to consider the symmetric version of the $M$-Wright function.
  \\ 
   Hereafter we provide
  a list of the main properties of this function, 
  which can be derived from the Laplace and Fourier transforms 
  for the corresponding Wright $M$-function
  in one variable.  
 \\
   From Eqs. (39) and (43) 
   we   derive the Laplace transform
    of $\MM_\nu(x,t)$ with respect to $t \in\RR^+$,   
   \begin{equation}
   \L\left\{\MM_\nu (x,t);  t \to s \right\}= s^{\nu-1}\, \e^{\ds \, -xs^\nu}\,.
\end{equation}  
    From Eq. (18) we   derive the Laplace transform of $\MM_\nu(x,t)$ with respect to $x\in \RR^+$,
\begin{equation}	
	\L\left\{\MM_\nu(x,t);  x \to s \right\}= E_{\nu}\left( -s t^\nu \right)\,.
\end{equation}	
     From Eq. (55) we   derive the Fourier transform
	 of $\MM_\nu(|x|,t)$ with respect to $x\in \RR$,
\begin{equation}	 
\F\left\{\MM_\nu(|x|,t);  x \to \kappa \right\}= 2E_{2\nu}\left( -\kappa^2 t^\nu \right)\,.
\end{equation}
 Using the  Mellin transforms,
   Mainardi et al. \cite{Mainardi-Pagnini-Gorenflo FCAA03}
derived 
the following interesting  integral formula of composition,
  \begin{equation}
 \MM_\nu(x,t)= \int_0^\infty \MM_\lambda(x,\tau)\, \MM_\mu(\tau,t)\, d\tau\,,\q \nu = \lambda \mu\,.
 \end{equation}
 Special cases of the  Wright $\MM$-function are simply derived for $\nu=1/2$ and $\nu=1/3$ from
 the corresponding ones in the complex domain, see Eqs. (28)-(29).
 We devote particular attention to the case $\nu=1/2$  for which we get  
 the Gaussian density in $\RR$,
\begin{equation} 
  \MM_{1/2}(|x|,t) = \rec{2\sqrt{\pi}t^{1/2}}\, \e^{\ds\, -x^2/(4t)} \,.
\end{equation}  
 For the limiting case $\nu=1$ we obtain
\begin{equation} 
 \MM_1 (|x|,t) = \rec{2} \left[\delta(x-t)+ \delta(x+t)\right]\,.
\end{equation} 
\section{The four sisters}
In this section we show how some  Wright functions of the second kind 
can provide an interesting generalization  of the three sisters discussed 
in Appendix A. 
The starting point is a  (not well- known) 
paper published  in 1970  by  Stankovic \cite{Stankovic BEOGRAD1970},
where  (in our notation) the following Laplace transform pair is proved 
rigorously:
\begin{equation}
 t ^{\mu-1}\, W_{{-\nu}, \mu}(x,t) 
\div  s^{-\mu}\, \e^{-xs^{\nu}}\, 
\q 0<\nu<1 \,, \q \mu \ge 0\,,
\label{eq:foursisterscompact}
\end{equation}
where $x$ and $t$ are positive.
We note that the Stankovic formula can be derived in a formal way by developing the 
exponential function in positive power of $s$ and inverting term by term
as described in the Appendix F of the book by Mainardi 
\cite{Mainardi BOOK2010}.
\\
We recognize  that the Laplace Transforms of the Three Sisters functions
 $\widetilde{\phi}(x,s)$, $\widetilde{\psi}(x,s)$ and  $\widetilde{\chi}(x,s)$ are particular cases of the
 Eq. (\ref{eq:foursisterscompact})  for $\nu= 1/2$, that is of 
 \begin{equation}
   t ^{\mu-1}\, W_{{-1/2}, \mu}(x,t)  \div    s^{-\mu}\, \e^{-x\sqrt{s}},
\label{eq:threesisterscompact}
\end{equation}
 according to the following scheme:
 \begin{itemize}
  \item[-]     $\widetilde{\phi} (x,s)$  with $\mu = 1$,
  \item[-]     $\widetilde{\psi} (x,s)$ 	with $\mu = 0$,
 \item[-]     $\widetilde{\chi} (x,s)$  with 	$\mu = 1/2 $.
 \end{itemize}
If $\nu$  is no longer restricted to $\nu = 1/2$  we define  
 \emph{Four Sisters functions}
 as follows 
\begin{equation}
\begin{split}
\mu = 0,~~~ \e^{-xs^{\nu}} &\div t^{-1}W_{-\nu, 0}(-xt^{-\nu}),\\
\mu = 1 - \nu, ~~~ \frac{\e^{-xs^{\nu}}}{s^{1-\nu}} &\div t^{-\nu}W_{-\nu, 1-\nu}(-xt^{-\nu}),\\
\mu = \nu, ~~~ \frac{\e^{-xs^{\nu}}}{s^{\nu}} &\div t^{\nu -1}W_{-\nu, \nu}(-xt^{-\nu}),\\
\mu = 1, ~~~ \frac{\e^{-xs^{\nu}}}{s} &\div W_{-\nu, 1}(-xt^{-\nu}). \\
\end{split} 
\end{equation}
\\
In the next figures  we show some plots of these functions, both in the $t$ and in the $x$ domain for some values of $\nu$ ($\nu = 1/4, 1/2, 3/4$).
\\
 Note that for $\nu = 1/2$ we only find three functions, (the Three Sisters functions)  of Appendix A
\begin{figure}[h!]
	\centering
	\includegraphics[width=105mm, height=45mm]{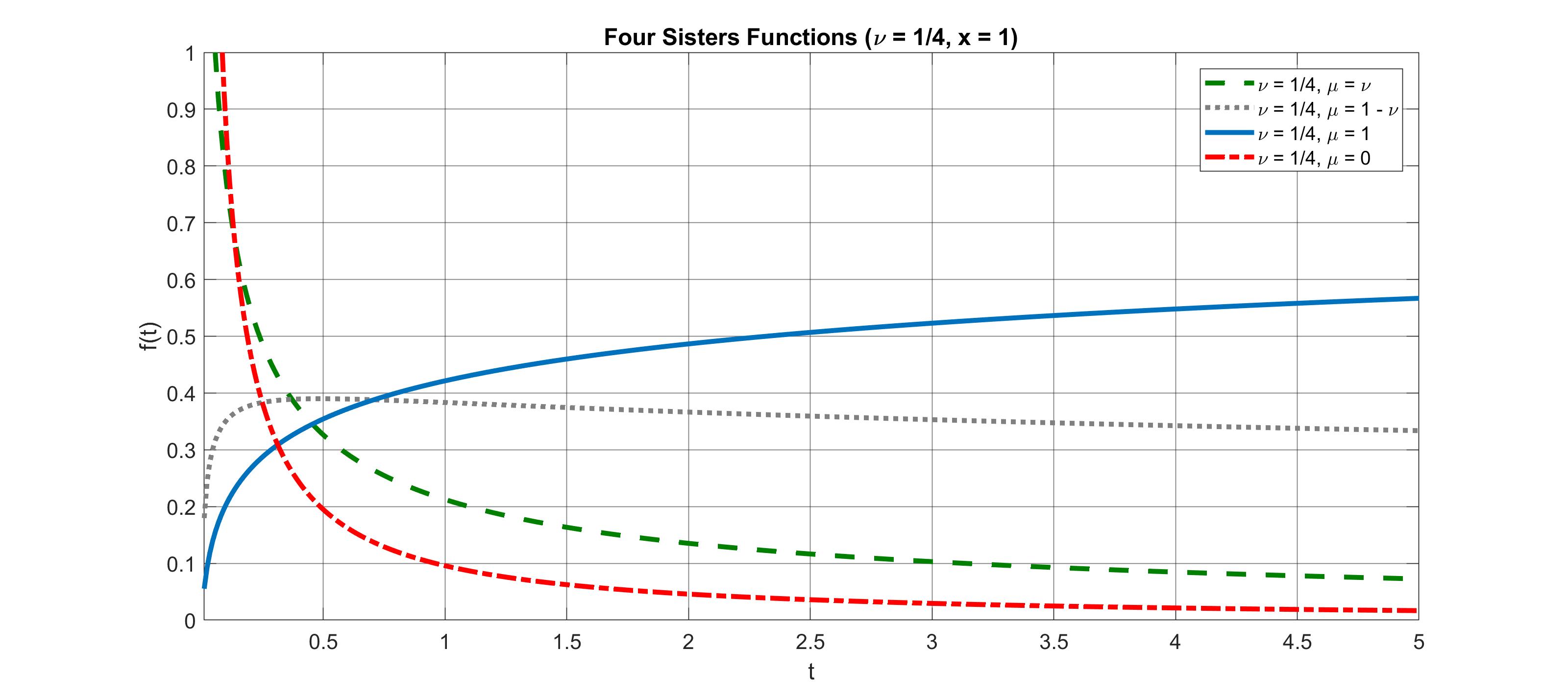}%
	\qquad\qquad
	\includegraphics[width=105mm, height=45mm]{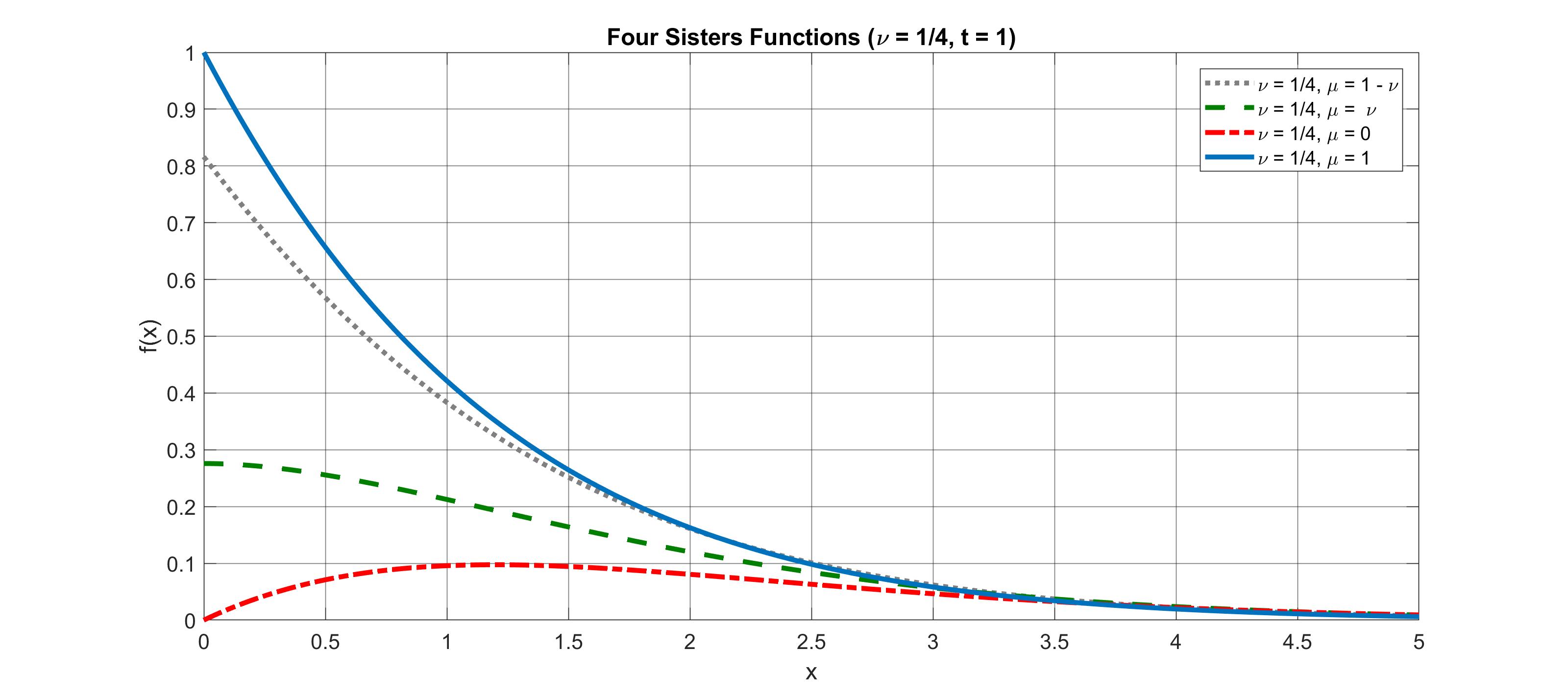}
	\caption{The four sisters functions in linear scale with $\nu = 1/4$; top: versus $t$ ($x=1$), bottom: versus $x$ ($t=1$)}
\end{figure}
\begin{figure}[h!]
	\centering
	\includegraphics[width=105mm, height=45mm]{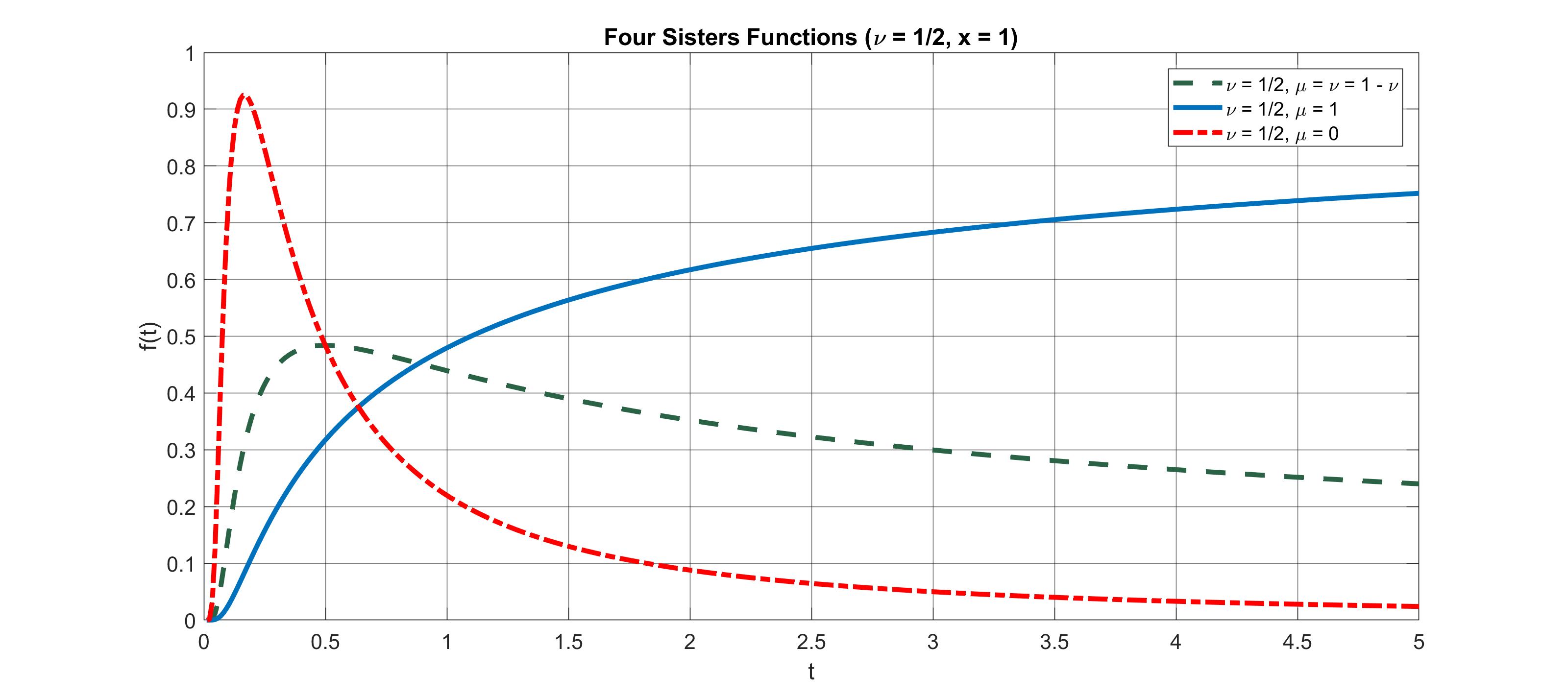}%
	\qquad\qquad
	\includegraphics[width=105mm, height=45mm]{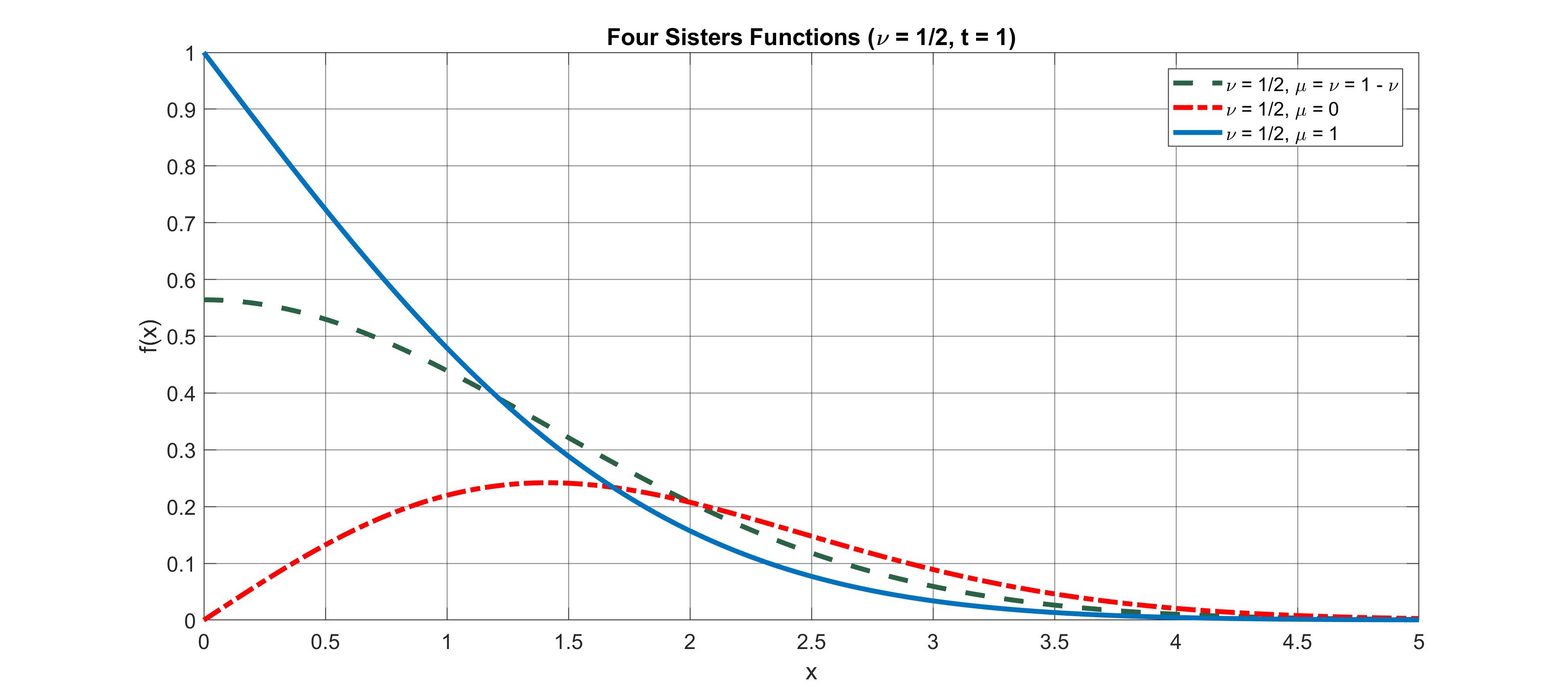}
	\caption{The three sisters functions in linear scale with $\nu = 1/2$; top: versus $t$ ($x=1$), bottom: versus $x$ ($t=1$)}
\end{figure}
\begin{figure}[h!]
	\centering
	\includegraphics[width=105mm, height=45mm]{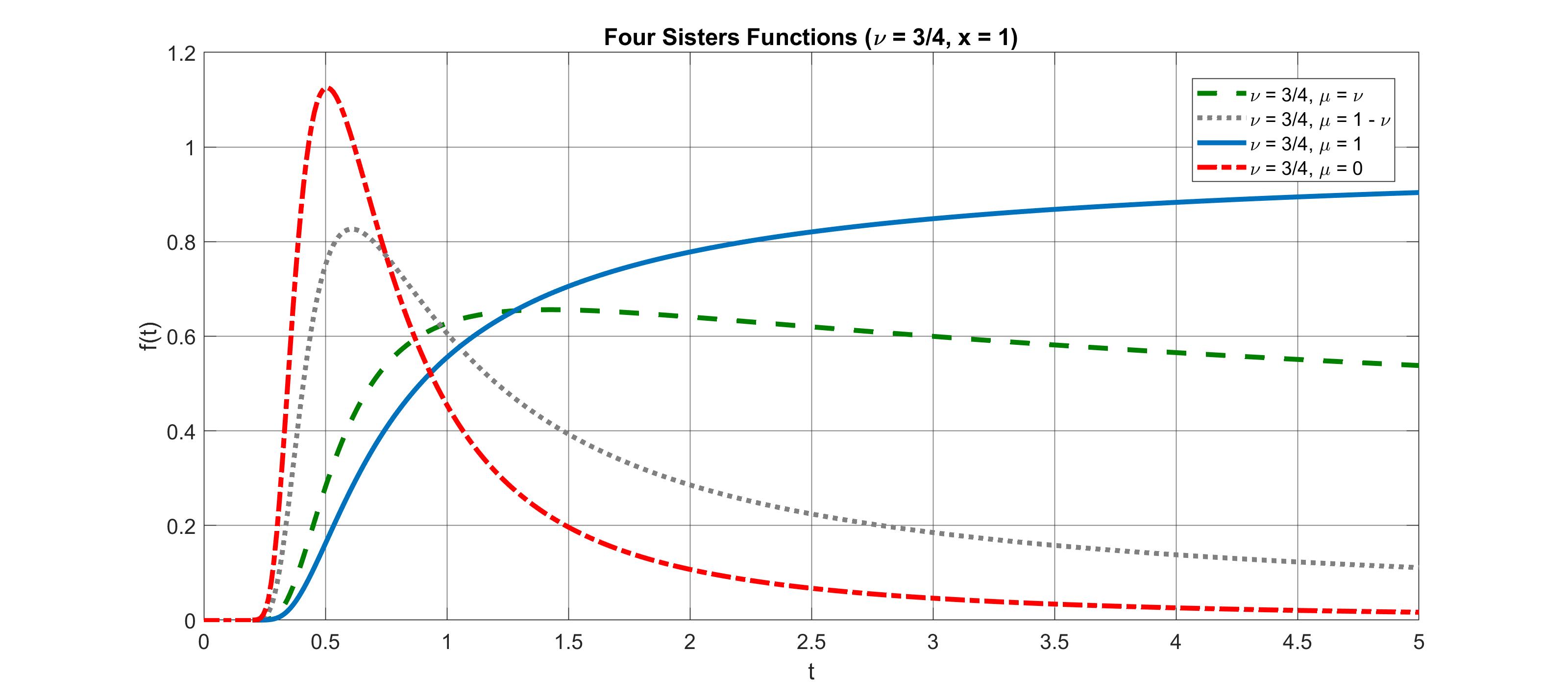}%
	\qquad\qquad
	\includegraphics[width=105mm, height=45mm]{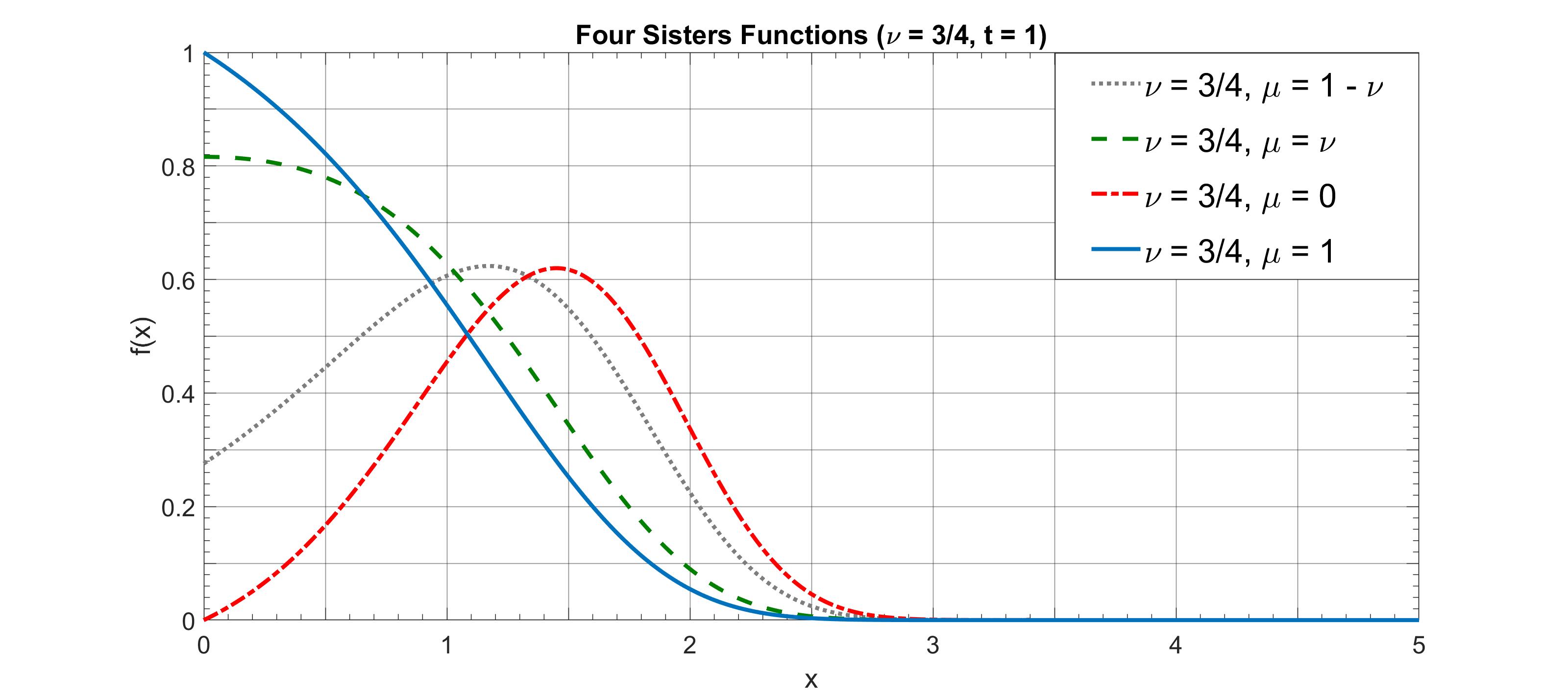}
	\caption{The four sisters functions in linear scale with $\nu = 3/4$; top: versus $t$ ($x=1$), bottom: versus $x$ ($t=1$)}
\end{figure}
\section{Conclusions}
\noindent

In our survey on the Wright functions we have distinguished   two kinds, 
pointing out  the  particular class  of the second kind.
Indeed these functions have been shown to play key roles in several processes 
 governed by non Gaussian processes, including sub-diffusion, transition to wave propagation, L{\'e}vy stable distributions.
 \\
Furthermore, we have devoted our attention to four functions of this class that we 
 agree to  called \emph{the Four Sisters functions}. 
 \\
All these items  justify  the relevance of the Wright functions  of the second kind in Mathematical Physics.
\\ 

\section*{Acknowledgments}{The research activity of both the authors
has been carried out in the framework of the activities of the National Group of Mathematical Physics (GNFM, INdAM).}

\section*{Appendix A: The standard diffusion equation and the three sisters}
In this Appendix let us  recall the Diffusion Equation 
in the one-dimensional case 
$$	\frac{\partial u}{\partial t} = D\frac{\partial^2 u}{\partial x^2}\,,
	\label{eq:diffusionequation}
\eqno(A.1)$$
where the constant   $D>0$ is the diffusion coefficient coefficient, whose dimensions are $L^2 T^{-1}$ and $x$ , $t$ denote the space and time coordinates, respectively.
\\ 
Two basic problems for Eq. (\ref{eq:diffusionequation}) are the \emph{Cauchy}
and \emph{Signalling} 
ones introduced hereafter
In these problems   some initial values and boundary conditions are set; specify the values attained by the field variable and/or by some of its derivatives on the boundary of the space-time domain is an essential step to guarantee the existence, the uniqueness and the determination of a solution of physical interest to the problem, not only for the Diffusion Equation.\\
\\
Two \emph{data functions} $f(x)$ and $g(t)$ are then introduced to write formally these conditions; some regularities are required to be satisfied by $f(x)$ and $g(t)$, and in particular $f(x)$ must admit the Fourier transform or the Fourier series expansion if the support is finite, while $h(t)$ must admit the Laplace Transform.\\ 
We also require without loss of generality that the field variable $u(x,t)$ is vanishing for $t < 0$ for every $x$ in the spatial domain.
Given these premises, we can specify the two aforementioned problems.\\
In the \emph{Cauchy problem} the medium is supposed to be unlimited
 ($-\infty < x < +\infty$) and to be subjected at $t = 0$ to a known disturbance provided by the data function $f(x)$. Formally:
$$
\begin{cases}
\lim_{t \rightarrow 0^+}{u(x, t)} = f(x), & -\infty < x < +\infty;\\
\lim_{x \rightarrow \pm \infty}{u(x, t)} = 0, & ~~~~~~~~~ t > 0.
\end{cases}  
\label{eq:cauchyproblem}     
\eqno(A.2)$$
This is a pure \emph{initial-value problem} (IVP) as the values are specified along the boundary $t=0$.
\\
In the \emph{Signalling problem} the medium is supposed to be semi-infinite
 ($0 \leq x < +\infty$) and initially undisturbed. At $x=0$  (the accessible end) and for $t > 0$ the medium is then subjected to a known disturbance provided by the causal function $g(t)$. Formally:
$$
\begin{cases}
\lim_{t \rightarrow 0^+}{u(x, t)} = 0, & 0 \leq x < +\infty;\\
\lim_{x \rightarrow 0^+}{u(x, t)} = g(t),~ \lim_{x \rightarrow  +\infty}{u(x, t)} = 0 & ~~~~~~~~~ t > 0.
\end{cases} 
\label{eq:signallingproblem}      
\eqno(A.3) $$
This problem is referred to as an \emph{initial boundary value problem} (IBVP) in the quadrant $\{x,t\}>0$.
\\
For each problem the solutions turn out to be expressed by a proper convolution
  between the data functions and the \emph{Green functions} $\mathcal{G}$, that are the fundamental solutions of the problems.
\\
For the Cauchy problem we have:
$$ 
	u(x, t) = \int_{-\infty}^{+\infty}{\mathcal{G}_C (\xi, t)f(x-\xi)d\xi} = \mathcal{G}_C (x, t)\ast f(x)\,,
		\label{eq:cauchyconvolution}
\eqno(A.4)$$   
with
$$ 
	\mathcal{G}_C (x, t) = \frac{1}{2\sqrt{\pi D t}}\e^{-x^2 /(4Dt)}.
	\label{eq:Green-Cauchy}
\eqno(A.5) $$ 
\\
For the Signalling problem we have:
$$ 
	u(x, t) = \int_{0}^{t}{\mathcal{G}_S (x, \tau)g(t-\tau)d\tau} = \mathcal{G}_S (x, t)\ast g(t)\,, \q -\infty<x<+\infty, \q t\ge0\,,
	\label{eq:signallingconvolution}
\eqno(A.6) $$ 
with
$$ 
	\mathcal{G}_S (x, t) = \frac{x}{2\sqrt{\pi D t^3}}\e^{-x^2 /(4Dt)}\,,
	\q x\ge 0, \q t\ge 0\,.
	\label{eq:Green-Signalling}
\eqno(A.7) $$ 
Following the lecture notes in Mathematical Physics by Mainardi
\cite{Mainardi DIFFUSION2019}, 
we note that the following relevant property is valid for $\{x,t\} >0$:
$$ 
	x \mathcal{G}_C(x,t) = t \mathcal{G}_S (x, t) = F(z)
	\label{eq:diffusionreciprocityrelation}
\eqno(A.8) $$ 
where
$$ 
	z = \frac{x}{\sqrt{Dt}}, ~~~~ F(z) = \frac{z}{2} M(z), ~~~~ M(z) = \frac{1}{\sqrt{\pi}} \e^{-z^2 /4}.
	\label{eq:reciprocityrelationstandarddiff}
\eqno(A.9)$$  
According to Mainardi' s notations, Eq. (\ref{eq:diffusionreciprocityrelation})
 is known as \emph{reciprocity relation}, $F(z)$ and $M(z)$ are called \emph{auxiliary functions} and $z$ is the \emph{similarity variable}.
\\
A particular case of the Signalling problem is obtained when $g(t) = H(t)$
(the Heaviside unit step function)   and the solution $u(x,t)$ turns out to be expressed in terms of the \emph{complementary error function}:
$$ 
	u(x,t) = \mathcal{H}_S (x, t) =
 \int_{0}^{t}{\mathcal{G}_S (x, \tau)d\tau}	=
	 \hbox{erfc} \Bigl( \frac{x}{2\sqrt{Dt}}\Bigr)\,, 
	 \q x\ge 0, \q t\ge 0\,.
\eqno(A.10) $$ 
\\
As well known,  the three  above fundamental solutions can be  obtained via the Fourier and Laplace  transform methods.
Introducing the parameter 
$a = |x|/\sqrt{D}$ 
the Laplace transforms of these  functions turns out to be simply 
related in the Laplace domain
$\mathrm{Re}(s) > 0$, as follows
$$ 
	\phi (a, t) := \hbox{erfc}\Bigl( \frac{a}{2\sqrt{t}}\Bigr) \div \frac{e^{-as^{1/2}}}{s} := \widetilde{ \phi} (a, s) ,
\label{eq:stepresponselaplacepair}
\eqno(A.11)$$ 
$$ 
	\psi (a, t) := \frac{a}{2\sqrt{\pi}}t^{-3/2}e^{-a^2 /(4t)} \div e^{-as^{1/2}} := \widetilde{ \psi}(a,s) , 
\label{eq:signallingproblemlaplacepair}
\eqno(A.12)$$ 
$$ 
\chi (a, t) := \frac{1}{\sqrt{\pi}}t^{-1/2}e^{-a^2 /(4t)} \div \frac{e^{-as^{1/2}}}{s^{1/2}} := \widetilde{ \chi}(a,s)
\label{eq:cauchyproblemlaplacepair}
\eqno(A.13)$$ 
where the sign $\div$ is used for the juxtaposition of a function with its Laplace transform.
We easily note that 
Eq. (\ref{eq:stepresponselaplacepair}) is related to the Step-Response problem, Eq. (\ref{eq:signallingproblemlaplacepair}) is related to the Signalling problem and Eq. (\ref{eq:cauchyproblemlaplacepair}) is related to the Cauchy problem.
Following the lecture notes by Mainardi 
\cite{Mainardi DIFFUSION2019} we agree to call the above functions
{\it the three sisters functions} for their role in the standard diffusion equation.
They will be discussed with details hereafter. 
\\
Everything that we have said above  will be found again as a special case of the \emph{Time Fractional Diffusion Equation} where the time derivative of the first order  is replaced by a suitable time derivative of non-integer order. 
\\
It is easy to demonstrate that each of them can be expressed as a function of one of the 2 others 
\textit{three sisters} (table \ref{threesisters-relations}).
\begin{table}[!htbp]
\renewcommand\arraystretch{2.3}
\begin{center}
\begin{tabular}{|p{2.1cm}|p{2.1cm}|p{2.1cm}|p{2.1cm}|}
\hline
   &   $ \widetilde{\phi} $   &   $ \widetilde{\psi} $  &  $ \widetilde{\chi} $  \\
\hline
$ \widetilde{\phi} $   &   $ \dfrac{\textrm{e}^{-a\sqrt{s}}}{s} $   &  $ \dfrac{\widetilde{\psi}}{s} $  & $ - \dfrac{1}{s} \dfrac{\partial\widetilde{\chi}}{\partial a} $ \\
\hline
$ \widetilde{\psi} $  &  $ s\; \widetilde{\phi} $  &  $ \textrm{e}^{-a\sqrt{s}} $   & $ -\dfrac{\partial\widetilde{\chi}}{\partial a} $   \\
\hline
$ \widetilde{\chi} $ &  $ -\dfrac{\partial\widetilde{\phi}}{\partial a} $   &   $ - \dfrac{2}{a}\dfrac{\partial\widetilde{\psi}}{\partial s} $  &  $ \dfrac{\textrm{e}^{-a\sqrt{s}}}{\sqrt{s}} $  \\
\hline
\end{tabular}
\end{center}
\caption{Relations among the \textit{three sisters} in the Laplace domain.}
\label{threesisters-relations}
\end{table}
\\
The \textit{three sisters} in the $ t $ domain may be all directly calculated by making use of the \textit{Bromwich formula} taking account of the contribution of the branch cut of $\sqrt {s}$ and of the pole of $1/s$. Wee obtain:
 $$\widetilde{\phi}(a,s) \div \phi (a,t) = 1-\dfrac{1}{\pi} \int_{0}^{\infty}
 \textrm{e}^{-rt} \sin (a\sqrt{r})\dfrac{\diff r}{r}\,,
 \eqno(A.11a)$$
 $$\widetilde{\psi}(a,s) \div \psi(a,t)=\dfrac{1}{\pi}\int_{0}^{\infty} \textrm{e}^{-rt} \sin (a\sqrt{r}) \diff r \,,
 \eqno(A.12a)$$
 $$\widetilde{\chi}(a,s)\div\chi(a,t)=\dfrac{1}{\pi}\int_{0}^{\infty} \textrm{e}^{-rt} \cos (a\sqrt{r}) \dfrac{\diff r}{\sqrt{r}} \,.
 \eqno(A.13a)$$
Then, through the substitution $ \rho=\sqrt{r} $, we arrive at the Gaussian integral and, consequently, we find the previous  explicit expressions of the \textit{three sisters}, that is:
 $$\phi (a,t) = \textrm{erfc}(\dfrac{a}{2\sqrt{t}}) = 1-\dfrac{2}{\sqrt{\pi}}\int_{0}^{a/2\sqrt{t}} \textrm{e}^{-u^{2}}\diff u\,,
 \eqno(A.11) $$   
 $$\psi(a,t)=\dfrac{a}{2\sqrt{\pi}} \: t^{-3/2} \: \textrm{e}^{-a^{2}/4t}\,,
 \eqno(A.12) $$
 $$ \chi(a,t)=\dfrac{1}{\sqrt{\pi}} \: t^{-1/2} \: \textrm{e}^{-a^{2}/4t} \;,
\eqno(A.13) $$
reminding the definition of the complementary error function.
 \\
 Alternatively, we can compute the \textit{three sisters} in  $ t $ domain by using the relations among the \textit{three sisters} in the Laplace domain listed in table \ref{threesisters-relations}. But in this case one of the \textit{three sisters} in 
  $ t $ domain must be already known. 
  Assuming to know $ \phi(a,t) $ 
from  Eq.  (A.11), we get:
\\  \\
- $ \psi(a,t) $ from $ \widetilde{\psi}(a,s) = s \: \widetilde{\phi}(a,s) $.
  Indeed,  noting
$$ s\;\widetilde{\phi}(a,s)\div \dfrac{\partial}{\partial t} \: \phi(a,t)$$
since $ \phi(a,0^{+})=0 $ 
we can obtain (A.12), namely
$$\psi(a,t)=\dfrac{a}{2\sqrt{\pi}} \: t^{-3/2} \;\textrm{e}^{-a^{2}/4t}\,;
$$
\\
- $ \chi(a,t) $ from $ \widetilde{\chi}(a,s) = -\dfrac{\partial}{\partial a}\widetilde{\phi}(a,s) $ where $ a $ is seen as  a parameter,. Indeed  it immediately follows (A.13), namely
$$
\chi(a,t)=-\dfrac{\partial}{\partial a} \;\phi(a,t)=\dfrac{1}{\sqrt{\pi}} \: t^{-1/2} \: \textrm{e}^{-a^{2}/4t}\,. $$
\\
For more details we refer the reader again to 
\cite{Mainardi DIFFUSION2019}.
 \section*{Appendix B: Essentials of Fractional Calculus}
Fractional calculus is the field of mathematical analysis which deals
with the investigation and applications of integrals and
derivatives of arbitrary order.
The term {\it fractional} is a misnomer, but it is
retained 	for historical reasons,
following the  prevailing use.
\vsp
This appendix is  based on the 1997 surveys by Gorenflo and Mainardi
\cite{Gorenflo-Mainardi CISM97} and by Mainardi \cite{Mainardi CISM97}.
 For more details on the classical treatment of	fractional
calculus
the reader is referred to 
the nice and rigorous book by Diethelm \cite{Diethelm LNM2010}
published in 2010 by Springer in the series Lecture Notes in Mathematics.
\vsp
According to the Riemann-Liouville approach  to fractional calculus,
the notion of fractional integral of order $\alpha$
($\alpha >0$)
 is a natural consequence
of the well known formula (usually attributed to Cauchy),
that reduces the calculation of the $n-$fold primitive of a function
$f(t)$ to a single integral of convolution type.
In our notation the Cauchy formula reads
$$
    J^n f(t) := f_{n}(t)=
 \rec{(n-1)!}\, \int_0^t \!\!  (t-\tau )^{n-1}\,f(\tau) \, d\tau\,,
    \q t > 0\,,\q n \in \NN   \,, \eqno(B.1) $$
where $\NN$ is the set of positive integers.
From this definition we note that $f_{n}(t)$
vanishes at $t=0$ with its
derivatives of order $1,2, \dots, n-1\,. $
For convention we require that	$f(t)$ and henceforth
$f_{n}(t)$ be a {\it causal} function, \ie identically
vanishing for $t<0\,. $
\vsp
In a natural way  one is  led
to extend the above formula
from positive integer values of the index to any positive real values
by using the Gamma function.
Indeed, noting that $(n-1)!= \Gamma(n)\,, $
and introducing the arbitrary {\it positive} real number
 $\alpha\,, $
one defines  the
 \ \underbar{{\it Fractional Integral of order} $\alpha >0 $} :
$$
J^\alpha \,f(t) :=
      \rec{\Gamma(\alpha )}\,
 \int_0^t (t-\tau )^{\alpha -1}\, f(\tau )\,d\tau \,,
   \q t > 0\,,\q \alpha  \in \RR^+
 \,,  \eqno(B.2) $$
where $\RR^+$ is the set of positive real numbers.
For complementation we define
$J^0 := I\, $ ({Identity operator)}, \ie we mean
$J^0\, f(t) = f(t)\,. $ Furthermore,
by $J^\alpha f(0^+)$ we mean the limit (if it exists)
of $J^\alpha f(t)$ for $t\to 0^+\,;$ this limit may be infinite.
\vsp
We note the {\it semigroup property}
$J^\alpha J^\beta = J^{\alpha +\beta}\,,
   \; \alpha\,,\;\beta	\ge 0\,,$
which implies the {\it commutative property}
$J^\beta  J^\alpha= J^\alpha J^\beta\,,$
and  the effect of our operators $J^\alpha$
on the power functions
$$
J^\alpha t^\gamma ={\Gamma (\gamma +1)\over \Gamma(\gamma +1+\alpha)}\,
		   t^{\gamma+\alpha}\,, \q \alpha \ge 0\,,
  \q \gamma >-1\,, \q t>0\,.
\eqno (B.3)
$$
These properties  are of course a natural generalization
of those known when the order is a positive integer.
\vsp
Introducing
the Laplace transform by the notation
$ {\mathcal{L}}\, \l\{  f(t) \r\}  := \int_0^\infty \!\!
   \e^{-st}\, f(t)\, dt = \widetilde f(s)\,, \; s \in \CC\,,$
and  using the sign $\div$ to denote a Laplace transform pair,
\ie
$ f(t) \div  \widetilde f(s) \,, $
we note the following rule for the   Laplace transform of
the fractional integral,
$$	   J^\alpha \,f(t) \div
     {\widetilde f(s)\over s^\alpha}\,,\q \alpha \ge 0\,,  \eqno(B.4)$$
which is the generalization
of the case with an $n$-fold repeated integral. 
\vsp
After the notion of fractional integral,
that of fractional derivative of order $\alpha$
($\alpha >0$)
becomes a natural requirement and one is attempted to
substitute $\alpha $ with $-\alpha $ in the above formulas.
However, this generalization  needs some care  in order to
guarantee the convergence of  the integrals   and
preserve the
well known properties of the ordinary derivative of integer
order.
\vsp
 Denoting by $D^n\,$ with $ n\in \NN\,, $
the operator of the derivative of order $n\,,$	we first note that
$ D^n \, J^n = I\,, \;	 J^n \, D^n \ne I\,,\q n\in \NN \,,
$
\ie $D^n$ is left-inverse (and not right-inverse) to
the corresponding integral operator $J^n\,. $
In fact we easily recognize from (B.1) that
$$  J^n \, D^n \, f(t) = f(t) - \sum_{k=0}^{n-1}
	f^{(k)}(0^+) \, {t^k\over k!}\,, \q t>0\,. \eqno(B.5)$$
As a consequence we expect that $D^\alpha $ is defined as left-inverse
to $J^\alpha $.  For this purpose, introducing the positive
integer $m$ such that $m-1 <\alpha \le m\,, $
one defines the
 \underbar{{\it Fractional Derivative of order} $\alpha >0 $}
as $\; D^\alpha \,f(t) := D^m \, J^{m-\alpha} \, f(t)\,,$
\ie
$$ \!\!\!
 D^\alpha \,f(t) :=
 \left\{
\begin{split}
 & {\ds {d^m\over dt^m}}\left[
  {\ds \rec{\Gamma(m-\alpha)}\int_0^t
    {f(\tau )\over (t-\tau )^{\alpha +1-m}} \,d\tau}\right] ,
 \qq  m-1 <\alpha < m, \\
   &  {\ds {d^m\over dt^m}} f(t)\,,
    \qq  \alpha =m\,.  
    \end{split}
    \right.
   \eqno(B.6) $$
Defining for complementation $D^0 = J^0 =I\,, $ then
we easily recognize that
$ D^\alpha \, J^\alpha = I \,,$  $\, \alpha \ge 0\,,$
and
$$ D^{\alpha}\, t^{\gamma}=
   {\Gamma(\gamma +1)\over\Gamma(\gamma +1-\alpha)}\,
     t^{\gamma-\alpha}\,,
 \q \alpha \ge 0\,,
  \q \gamma >-1\,, \q t>0\,.
\eqno (B.7)
$$
Of course, these properties are a natural generalization
of those known when the order is a positive integer.
\vsp
Note the remarkable fact that the fractional derivative $D^\alpha\, f$
is not	zero
for the constant function $f(t)\equiv 1$ if $\alpha \not \in {\NN}\,. $
In fact, (B.7) with $\gamma =0$ teaches us that
$$
D^\alpha 1 = {t^{-\alpha}\over \Gamma(1-\alpha)}\,,\q \alpha\ge 0\,,
\q t>0\,.  \eqno (B.8)
$$
This, of course, is $\equiv 0$ for $\alpha \in{\NN}$, due to the
poles of the gamma function in the points $0,-1,-2,\dots$.
We now observe that an alternative definition
of fractional derivative was  introduced by Caputo 
in 1967 \cite{Caputo 1967} in a geophysical journal
and in 1969 \cite{Caputo BOOK1969} in a book in Italian.
Then the Caputo  definition was adopted in 1971 by Caputo and Mainardi 
\cite{Caputo-Mainardi 1971a,Caputo-Mainardi 1971b}
in the framework  of the theory of {\it Linear Viscoelasticity}.
Nowadays it is usually referred  to as the {\it Caputo fractional derivative}
and reads 
$\; D_*^\alpha	\, f(t) :=   J^{m-\alpha}\, D^{m} \, f(t) $
with $m-1 <\alpha \le m\,, \; m\in \NN\,,$ \ie
$$
 D_*^\alpha \,f(t) :=
 \left\{
\begin{split}
  & {\ds \rec{\Gamma(m-\alpha)}}\,{\ds\int_0^t
 {\ds {f^{(m)}(\tau)\over (t-\tau )^{\alpha +1-m}}} \,d\tau} \,,
   \qq  m-1<\alpha <m\,, \\
    & {\ds {d^m\over dt^m}} f(t)\,,
     \qq  \alpha =m\,.
    \end{split}
    \right . 
   \eqno(B.9) $$
We note that there are a number of discussions on the priority of this 
definition that surely was formerly considered by Liouville as sated
by    Butzer and Westphal \cite{Butzer-Westphal 00}.
However Liouville did not recognize the relevance of this  
representation derived by a trivial integration by part whereas Caputo,
even if unaware of the Riemann-Liouville representation, promoted his 
definition in several papers over all for the applications where the  Laplace transform plays a fundamental role.
We agree to
denote (B.9) as the {\it Caputo fractional derivative}
to distinguish it from the standard Riemann-Liouville fractional
derivative (B.6). 
     \\
   The Caputo  definition (B.9) is of course more restrictive than the Riemann-Liouville definition   (B.6), in that
requires the absolute integrability of the  derivative of order $m$.
Whenever we use the operator   $D_*^\alpha$   we (tacitly) assume that
     this condition is met.
 We  easily recognize that in general
$$  D^\alpha\, f(t) := D^{m} \, J^{m-\alpha} \, f(t)
 \ne J^{m-\alpha}\, D^{m} \, f(t):= D_*^\alpha \, f(t)\,,
 \eqno(B.10)
 $$
 unless   the function	$f(t)$ along with its first $m-1$ derivatives
 vanishes at $t=0^+$.
In fact, assuming that
the passage of the $m$-derivative under
the integral is legitimate, one     
recognizes that,  for $ m-1 <\alpha  < m \,$  and $t>0\,, $
$$
    D^\alpha \, f(t) =
  D_*^\alpha   \, f(t) +
  \sum_{k=0}^{m-1}   {t^{k-\alpha}\over\Gamma(k-\alpha +1)}
    \, f^{(k)}(0^+) \,, \eqno(B.11)    $$
 and therefore, recalling the fractional derivative of the power
functions (B.7),
$$
   D^\alpha \left( f(t) -
 \sum_{k=0}^{m-1} {t^k \over k!} \, f^{(k)} (0^+)\right)
     =	D_*^\alpha  \, f(t)  \,.\eqno(B.12)  $$
The alternative definition (B.9) for the
fractional derivative  thus incorporates the initial values
of the function and of its integer derivatives of lower order.
The subtraction of the Taylor polynomial of degree $m-1$ at $t=0^+$
from $f(t)$ means  a sort of
regularization	of the Riemann-Liouville fractional derivative.
In particular for $0<\alpha<1$
we get
$$ 
   D^\alpha \left( f(t) -  f (0^+)\right)
     =	D_*^\alpha  \, f(t)  \,.  $$
 According to the Caputo  definition,
the relevant property for which the fractional derivative
of a constant is still zero can be easily recognized,
 \ie
$$ D_*^\alpha  1 \equiv 0\,,\q	 \alpha >0\,.\eqno(B.13)$$
   \vsp
We now explore the most relevant differences between the two
fractional derivatives (B.6) and (B.9). 
We observe, again by looking at (B.7), that
$D^\alpha t^{\alpha -1} \equiv 0\,, \; \alpha>0\,, \; t>0\,.$
From above  we thus recognize
the following statements about functions
which  for $t>0\, $   admit the same fractional derivative
of    order $\alpha \,, $
with $m-1 <\alpha \le m\,,$ $\; m \in \NN\,, $
$$    D^\alpha \, f(t) = D^\alpha  \, g(t)
   \,  \Longleftrightarrow  \,
  f(t) = g(t) + \sum_{j=1}^m c_j\, t^{\alpha-j} \,,
    \eqno(B.14) $$
$$    D_* ^\alpha \, f(t) = D_*^\alpha	\, g(t)
   \,  \Longleftrightarrow  \,
  f(t) = g(t) +  \sum_{j=1}^m c_j\, t^{m-j} \,.
    \eqno(B.15) $$
In these formulas the coefficients $c_j$ are arbitrary constants.
\vsp
For the two definitions we also note a difference
with respect to the {\it formal} $\,$ limit  as
 $\alpha \to {(m-1)}^+$. From (B.6) and (B.9) we obtain
respectively,
$$ \alpha \to (m-1)^{+}\,\Longrightarrow\,
 D^\alpha \,f(t) \to	D^m\,  J\, f(t) = D^{m-1}\, f(t)
   \,; \eqno(B.16) $$
$$\alpha \to {(m-1)}^{+} \,\Longrightarrow\,
 D_*^\alpha \, f(t) \to J\, D^m\, f(t) =
       D^{m-1}\, f(t) - f^{(m-1)} (0^+)\,. \eqno(B.17) $$
\vsp
We now consider the {\it Laplace transform} of the two fractional
derivatives.
For the standard fractional derivative $D^\alpha $
the Laplace transform,	assumed to exist,  requires the knowledge of the
(bounded) initial values of the fractional integral $J^{m-\alpha }$
and of its integer  derivatives of order $k =1,2, \dots, m-1\,. $
The corresponding rule reads, in our notation,
$$ D^\alpha \, f(t) \div
      s^\alpha\,  \widetilde f(s)
   -\sum_{k=0}^{m-1}  D^k\, J^{(m-\alpha)}\,f(0^+) \, s^{m -1-k}\,,
  \q m-1<\alpha \le m \,. \eqno(B.18)$$
\vsp
The {\it Caputo fractional derivative} appears more suitable to
be treated by the Laplace transform technique in that it requires
the knowledge of the (bounded)
initial values of the function
and of its integer  derivatives of order $k =1,2, \dots, m-1\,, $
in analogy with the case when $\alpha =m\,. $
In fact,
by using (B.4) and noting that
$$\! \!\!\! J^\alpha  \, D_*^\alpha \, f(t) =
    J^\alpha\, J^{m-\alpha }\, D^m \, f(t) =
     J^m\, D^m \, f(t) = f(t) -
  \sum_{k=0}^{m-1} {f^{(k)}(0^+)} {t^k \over k!}
  . \eqno(B.19)$$
we easily prove  the following rule for the Laplace transform,
$$ D_*^\alpha \, f(t) \div
      s^\alpha\,  \widetilde f(s)
   -\sum_{k=0}^{m-1}  f^{(k)}(0^+) \, s^{\alpha -1-k}\,,
  \q m-1<\alpha \le m \,. \eqno(B.20)$$
Indeed, the  result (B.20), first stated by Caputo by using the
Fubini-Tonelli theorem, appears  as the most "natural"
generalization of the corresponding result well known for $\alpha =m\,. $
\\
In particular Gorenflo and Mainardi
have pointed out the major utility of the
Caputo fractional derivative
in the treatment of differential equations of fractional
order for {\it physical applications}.
In fact, in physical problems,	the initial conditions are usually
expressed in terms of a given number of bounded values assumed by the
field variable and its derivatives of integer order,
no matter if
the governing evolution equation may be a generic integro-differential
equation and therefore, in particular,	a  fractional differential
equation.
\section*{Appendix C: The L\'evy stable distributions}
\begin{figure}[ht!]
\centering
\includegraphics[width=.50\textwidth]{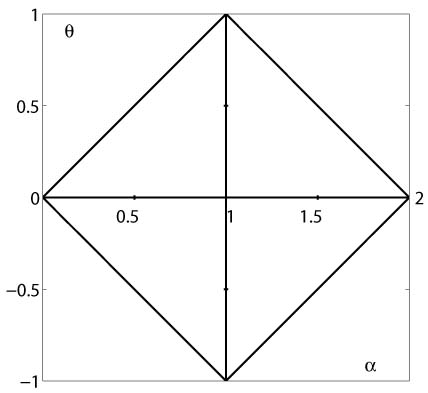}
 \caption{The Feller-Takayasu diamond for L\'evy stable densities. 
 \label{fig:F.4} }
\end{figure}

 \noindent
 We now introduce
 the so-called  {\it L\'evy stable  distributions}.. 
 The term stable has been assigned by the French  mathematician Paul L\'evy,
 who, in the tuenties of the last century, started a systematic research 
in order to generalize the celebrated {\it Central Limit Theorem}\index{Central limit theorem} to 
probability distributions  with infinite variance.  
For stable distributions we can assume the following 
 {\sc Definition}:
{\it If two independent real random variables
with the same shape or type of distribution are combined linearly and
the distribution of the resulting random variable has  the same shape,
the common distribution (or its type, more precisely) is said to be
stable}.
\vsp
The restrictive condition of stability enabled L\'evy (and then other authors) to derive
the {\it canonic form} for the characteristic function of the densities of these distributions.  
Here we follow the parameterization by Feller \cite{Feller 52,Feller BOOK71}
 revisited by Gorenflo \& Mainardi in \cite{GorMai FCAA98}, see also    \cite{Mainardi LUMAPA01}.
 Denoting by  $L_\alpha^\theta(x) $ a generic stable density  in $\RR$, 
 where $\alpha$ is the {\it index of stability} and  
 and $\theta$ the asymmetry parameter, improperly called {\it skewness}, 
 its characteristic function reads: 
$$
{\ds L_\alpha^\theta(x)}  \div {\ds \widehat{L}_\alpha ^\theta(\kappa)}  = 
  {\ds \exp \left[- \psi_\alpha ^\theta(\kappa )\right]} \,, \q
   {\ds \psi_\alpha  ^\theta(\kappa )} =
   {\ds |\kappa|^{\ds \alpha } \, \e^{\ds  i (\sgn \kappa)\theta\pi/2}}\,, 
 \eqno (C.1) $$
$$ 0<\alpha  \le 2\,, \;
 |\theta| \le  \,\hbox{min}\, \{\alpha  ,2-\alpha  \}\,.$$
 \vsp
 We note that the allowed region for the 
parameters $\alpha  $ and $\theta$
turns out to be
 a {diamond} in the plane $\{\alpha  , \theta\}$
with vertices in the points
$(0,0)\,, \,(1,1)\,, \, (1,-1) \,,\,(2,0)$,
that we call the {\it Feller-Takayasu diamond},
see Figure \ref{fig:F.4}. 
For values of $\theta$ on the border of the diamond
(that is $\theta = \pm \alpha $ if $0<\alpha  < 1$, and $\theta = \pm (2-\alpha )$ if $1<\alpha  <2$)
we obtain the so-called {\it  extremal stable densities}.
  \vsp 
\vsp
 We note the {\it symmetry relation}
$L_\alpha ^\theta (-x)=	L_\alpha ^{-\theta} (x)$, so that a stable density with $\theta=0$ is symmetric. 
  \vsp
  Stable distributions have noteworthy properties of which the interested reader 
  can be informed from the relevant existing literature. 
 Here-after we recall some  peculiar {\sc Properties}:
 \pni
- {\it The class of  stable distributions possesses its own {\it domain of attraction}}, see \eg \cite{Feller BOOK71}.
\pni
- {\it Any stable  density  is {unimodal} and indeed {bell-shaped}}, \ie
its $n$-th derivative has exactly $n$ zeros in $\RR$, 
see Gawronski \cite{Gawronski 84},
 Simon \cite{Simon 2015}
 and Kwa{\'s}nicki  \cite{Kwasnicki 2017}.
 \pni  
- {\it The stable distributions are  {self-similar} and  {infinitely divisible}}.
\vsp
These properties derive from the canonic form (C.1) through the scaling property of
the Fourier transform.
\pni 
{\it Self-similarity}  means
$$L_\alpha ^\theta (x,t) \div \exp \left [-t \psi _\alpha  ^\theta(\kappa)\right]
\Longleftrightarrow  L_\alpha ^\theta (x,t)  = t^{-1/\alpha }\,L_\alpha ^\theta (x/t^{1/\alpha} ) ]\,, 
 \eqno (C.2)  $$
where $t$ is a positive parameter.
If $t$ is time, then  $L_\alpha ^\theta (x,t)$ is a spatial density evolving on time with self-similarity. 
 \pni 
   {\it Infinite divisibility}  means that  for every positive integer
$n$,  the characteristic function  can be expressed as the $n$th power of some characteristic function, 
so that
any stable distribution can be expressed as the	$n$-fold convolution of a
stable distribution of the same type. 
Indeed, taking in (C.1) $\theta=0$, without loss of generality, we have
$$\e^{-t|\kappa|^\alpha } = \left[\e^{-(t/n)|\kappa|^\alpha }\right]^n 
\Longleftrightarrow  L_\alpha ^0 (x,t)  = 
\left[L_\alpha ^0 (x,t/n)\right]^{*n} \,, \eqno(C.3)$$  
where   
$$\left[L_\alpha ^0 (x,t/n)\right]^{*n} :=
L_\alpha ^0 (x,t/n) * L_\alpha ^0 (x,t/n) *  \dots * L_\alpha ^0 (x,t/n) $$
is the multiple Fourier convolution in $\RR$ with $n$ 
identical terms. 
\vsp
Only for a few particular cases, the inversion of the Fourier transform in (C.1)
can be carried out using standard tables, and well-known probability distributions are obtained.
\vsp
For $\alpha  =2$ (so $\theta =0$), we recover the  {\it Gaussian pdf}, that turns out to be the 
only stable density with finite variance, and more generally with finite  moments of any order 
$\delta \ge 0$. In fact
$$ L_2^0(x) = \frac{1}{2\sqrt{\pi}}\e^{\,\ds -x^2/4} 
\,.\eqno(C.4)$$
All the other stable densities have finite absolute moments  of order 
$\delta \in [-1, \alpha )$
as we will later show.
\vsp
For $\alpha  =1 $  and  $|\theta| <1$, we get 
$$L_1^\theta (x) =
 \rec{\pi} \, {  \cos (\theta \pi/2)
 \over [x+   \sin (\theta \pi/2)]^2 +[ \cos (\theta \pi/2)]^2 }\,,\eqno(C.5)$$
which  for $\theta=0$ includes the  {\it Cauchy-Lorentz pdf}.
$$ L_1^0(x) = \frac{1}{\pi} \frac{1}{1+x^2}   
\,.\eqno(C.6)$$
In the limiting cases   $\theta = \pm 1$ for $\alpha =1$ we obtain
 the {\it  singular Dirac pdf's} 
 $$ L_1^{\pm 1}(x)=\delta(x \pm 1)\,.\eqno(C.7)$$
\vsp
In general, we must recall the power series expansions 
provided in \cite{Feller BOOK71}.   
We restrict our attention to $x>0$
since the evaluations for $x<0$  can be obtained using the symmetry relation.
%
The convergent expansions of $L_\alpha ^{\theta} (x)$  ($x>0$) turn out to be;
 \\ 
for $ 0<\alpha   <1\,,\q |\theta| \le \alpha   \,:$
 $$L_\alpha  ^\theta (x) =
{1\over \pi\,x}\,  \sum_{n=1}^{\infty}
   (-x^{-\alpha  })^n \, {\Gamma (1+ n\alpha  )\over n!}\,
  \sin \left[{ n\pi\over 2}(\theta -\alpha  )\right]\,;
    \eqno(C,8)  $$
\\
for $ 1<\alpha   \le 2\,, \q |\theta| \le 2-\alpha  \,:$
$$ L_\alpha  ^\theta (x)=
{1\over \pi\,x}\,  \sum_{n=1}^{\infty}
   (-x)^{n} \, {\Gamma (1+ n/\alpha  )\over n!}\,
  \sin \left[{ n\pi\over 2\alpha   }(\theta -\alpha  )\right]\,.
 \eqno(C.9) $$
From the series in (C.8) and the  symmetry relation
we note that {\it the
extremal stable densities for $0<\alpha   <1$ are
unilateral}, precisely vanishing for $x>0$ if $\theta =\alpha $,
vanishing for $x<0$ if $\theta =-\alpha $.
In particular the unilateral extremal densities $L_\alpha ^{-\alpha }(x)$
 with $0<\alpha <1$ 
have support in $\RR^+$ and Laplace transform
$\exp (-s^\alpha )$. For $\alpha=1/2$ we obtain the so-called {\it L\'evy-Smirnov} $pdf$:
$$ L_{1/2}^{-1/2} (x) =
  {\,x^{-3/2}\over 2\sqrt{\pi}}\, \,\e^{\ds\, - 1/(4x)}\,,
 \q x \ge 0	\,. \eqno(C.10)$$
As a consequence of the convergence of the series
in (C.8)-(C.9) and of the symmetry relation
we recognize that
the stable $pdf$'s with $1< \alpha \le 2$
are  entire functions, whereas
 with $0< \alpha <1$
have the form
$$L_\alpha^\theta (x) = \begin{cases}
   (1/x) \,\Phi_1(x^{-\alpha })
	& \hbox{for} \; x>0\,, \\
  (1/|x|) \,\Phi_2(|x|^{-\alpha })
	& \hbox{for} \; x<0\,,
	\end{cases}
	  \eqno(C.11)		 $$
where $\Phi_1(z)$ and $\Phi_2(z)$ are distinct {entire functions}.
The  case $\alpha =1$ ($|\theta|< 1$) must be considered  in the limit for
$\alpha \to 1$ of  (C.8)-(C.9), because the corresponding
series reduce to power series akin with geometric series
in $1/x$ and $x$, respectively, with a finite radius of convergence.
The corresponding stable $pdf$'s are no longer represented by
entire functions, as can be noted directly from their explicit expressions 
(C.5)-(C.6).
\vsp
We omit to provide the asymptotic representations of the stable densities referring the interested reader to Mainardi et al (2001) 
\cite{Mainardi LUMAPA01}. However, based on asymptotic representations, we  can state as follows;
 for $0<\alpha <2$ the stable $pdf$'s  exhibit
{\it fat tails} in such a way that their  absolute moment
of  order $\delta$ is finite only if $-1 < \delta <\alpha $.
More precisely, one can show that for non-Gaussian,
not extremal, stable densities
the asymptotic decay of the tails  is
$$ L_\alpha^\theta (x )= O\left(|x|^{-(\alpha +1)}\right)\,, \q
	       x \to \pm \infty\,. \eqno(C.12)$$
For the extremal densities with $\alpha \ne 1$
this is valid only for one
tail (as $|x|\to \infty$), the other (as $|x|\to \infty$) being of exponential order.
For $1<\alpha <2$ the extremal $pdf$'s are two-sided and exhibit
an exponential left tail  (as $x \to -\infty)$
if $\theta  =+(2-\alpha)\,,  $
or  an exponential right tail  (as $x \to +\infty $)
if $\theta  =-(2-\alpha)\,.$
Consequently, the Gaussian $pdf$  is the unique
 stable density with finite variance.
Furthermore,  when $0<\alpha \le 1$,
the first absolute moment  is infinite  so
we should use the  median instead of the non-existent expected value
in order to characterize the corresponding $pdf$.
\vsp
Let us also recall a relevant identity between stable densities\index{L\'evy stable distribution}
 with index $\alpha $ and $1/\alpha$   (a sort of reciprocity relation) pointed out in \cite{Feller BOOK71},
that is, assuming $x>0$,
$$ \rec{x^{\alpha +1}}\, L_{1/\alpha}^\theta (x^{-\alpha} )
  =L_\alpha^{\theta^*} (x)\,,  \;1/2\le \alpha\le 1\,,\;
  \theta ^*=\alpha(\theta +1)-1 \,. \eqno(C.13)$$
 The condition $1/2\le \alpha \le 1$ implies $1\le  1/\alpha \le 2$. A check shows  
 that $\theta^*$ falls within the prescribed range
$|\theta ^*|\le\alpha$ if $|\theta |\le 2-1/\alpha $.
\\
We leave as an exercise for the interested reader the verification of this reciprocity relation 
in the limiting cases $\alpha=1/2$ and $\alpha=1$. 



\end{document}